\documentclass{siamart0516}
\usepackage{latexsym,amssymb,enumerate,amsmath} 
  \usepackage{graphicx} 
  \usepackage{tikz}
  \usepackage[square,sort,comma,numbers]{natbib}
  \usepackage{hyperref}
\usepackage{pgfplots}

\numberwithin{equation}{section}

  \def\black{\color{black}}
\newcommand{\dd}{\text{d}}
\def\eps{\varepsilon}

\def\P{\mathbb{P}}
\def\R{\mathbb{R}}
\def\L{\mathcal{L}}
\def\LL{\mathbb{L}}
\def\O{\mathcal{O}}
\def\RR{\mathcal{R}}
\def\x{\mathbf{x}}

\def\X{\mathbf{X}}

\def\Z{\mathbf{Z}}

\def\f{A}
\def\m{B}

\def\Dtr{D^{\textup{tr}}}
\def\Drot{D^{\textup{rot}}}
\def\Dsurf{D^{\textup{surf}}}
\def\Deff{D^{\textup{eff}}}
\def\ks{k_{\textup{smol}}}
\def\kg{k_{\textup{geo}}}
\def\ke{k_{\textup{emp}}}

\def\kb{k_{\textup{bind}}}
\def\keff{k_{\textup{eff}}}
\def\kbp{k_{\textup{bp}}}
\def\chiqc{\chi_{\textup{qc}}}
\def\chib{\chi_{\textup{b}}}

\usepackage{amsfonts}
\usepackage{graphicx}
\usepackage{epstopdf}
\usepackage{algorithmic}
\ifpdf
  \DeclareGraphicsExtensions{.eps,.pdf,.png,.jpg}
\else
  \DeclareGraphicsExtensions{.eps}
\fi

\newcommand{\TheTitle}{Bimolecular binding rates for pairs of spherical molecules with small binding sites}
\newcommand{\ShortTitle}{Bimolecular binding rates}
\newcommand{\TheAuthors}{Claire E. Plunkett and Sean D. Lawley}

\headers{\ShortTitle}{\TheAuthors}

\title{{\TheTitle}\thanks{
\funding{This work was supported by the National Science Foundation (DMS-1944574, DMS-1814832, and DMS-1148230).}}}

\author{Claire E. Plunkett\thanks{Department of Mathematics, University of Utah, Salt Lake City, UT 84112 USA (\texttt{plunkett@math.utah.edu}).}\and {Sean D. Lawley\thanks{Department of Mathematics, University of Utah, Salt Lake City, UT 84112 USA (\texttt{lawley@math.utah.edu}).}}}
\date{\today}

\ifpdf
\hypersetup{
  pdftitle={\TheTitle},
  pdfauthor={\TheAuthors}
}
\fi

\begin{document}

\maketitle


\begin{abstract}
Bimolecular binding rate constants are often used to describe the association of large molecules, such as proteins. In this paper, we analyze a model for such binding rates that includes the fact that pairs of molecules can bind only in certain orientations. The model considers two spherical molecules, each with an arbitrary number of small binding sites on their surface, and the two molecules bind if and only if their binding sites come into contact (such molecules are often called ``patchy particles'' in the biochemistry literature). The molecules undergo translational and rotational diffusion, and the binding sites are allowed to diffuse on their surfaces. Mathematically, the model takes the form of a high-dimensional, anisotropic diffusion equation with mixed boundary conditions. We apply matched asymptotic analysis to derive the bimolecular binding rate in the limit of small, well-separated binding sites. The resulting binding rate formula involves a factor that depends on the electrostatic capacitance of a certain four-dimensional region embedded in five dimensions. We compute this factor numerically by modifying a recent kinetic Monte Carlo algorithm. We then apply a quasi chemical formalism to obtain a simple analytical approximation for this factor and find a binding rate formula that includes the effects of binding site competition/saturation. We verify our results by numerical simulation. 
\end{abstract}



\section{Introduction}

The association of molecules to form dimers or larger complexes is characterized by bimolecular binding rate constants. To illustrate, consider two proteins, $A$ and $B$, which bind to form a complex, $C$. If $[A]$, $[B]$, and $[C]$ denote their respective concentrations, then the law of mass action \cite{keener09} implies that the concentration of the complex satisfies the ordinary differential equation (ODE),
\begin{align*}
\frac{\dd}{\dd t}[C]
=k [A] [B],
\end{align*}
for some bimolecular binding rate constant $k>0$ (also called a second-order rate constant). How does one determine $k$?

Protein-protein binding occurs through interactions between localized \emph{binding sites} on each protein. Hence, two commonly assumed \cite{zhou2010} conditions for protein-protein binding are a \emph{proximity} condition and an \emph{orientation} condition:
\begin{enumerate}[(i)]
\item
Proteins must be in sufficient proximity to bind.
\item
Binding sites must be properly oriented to bind.
\end{enumerate}
Smoluchowski's classical theory \cite{smoluchowski1917} provides a formula for the binding rate constant $k$ if we ignore the orientation condition (ii) (and this theory has had an indelible effect on how we understand binding kinetics \cite{gudowska2017, holcman2017}).

This classical theory involves the probability, $p(r)$, that two spherical proteins ($A$ and $B$) diffusing in three dimensions (3D) never bind to each other, given that they are initially separated by distance $r$. This probability satisfies Laplace's equation, 
\begin{align}\label{pdeG}
0=\frac{2}{r}\partial_{r}p+\partial_{rr}p,\quad \text{for }r>R,
\end{align}
where
\begin{align*}
R=R_{A}+R_{B}>0
\end{align*}
is the sum of the protein radii. In particular, $R$ is called the reaction radius and is the proximity in condition (i) at which the proteins bind. Since proteins that start far from each other will never bind, we obtain the far-field condition,
\begin{align}\label{ffc}
\lim_{r\to\infty}p=1.
\end{align}
The classical theory assumes that proteins bind immediately upon contact, which yields an absorbing boundary condition at the reaction radius,
\begin{align}\label{acG}
p=0,\quad \text{for }r=R.
\end{align}
The solution to \eqref{pdeG}-\eqref{acG} is simply $p(r)=1-R/r$. Calculating the flux at the reaction radius yields the classical Smoluchowski bimolecular binding rate constant, $k=\ks$,
\begin{align}\label{ksmol}
\ks
:=\Dtr\int_{r=R}\partial_{r}p\,\dd S
=4\pi \Dtr R.
\end{align}
where
\begin{align*}
\Dtr=\Dtr_{A}+\Dtr_{B}>0
\end{align*}
is the sum of the protein translational diffusivities. {\black Plugging typical values for proteins of $R\approx4\,\text{nm}$ and $\Dtr\approx2.5\times10^{8}\,\text{nm}^{2}\,\text{sec}^{-1}$ into \eqref{ksmol} yields that the Smoluchowski rate constant is on the order of \cite{northrup1992, zhou2010, zhou2013}}
\begin{align*}
\textcolor{black}{
\ks\approx7\times10^{9}\;\text{M}^{-1}\text{sec}^{-1}.
}
\end{align*}

Since this classical calculation ignores the orientation condition (ii) above, $\ks$ is an upper bound for binding rates \cite{zhou2010,zhou2013}. Indeed, $\ks$ tends to overestimate experimentally measured rates by several orders of magnitude \cite{northrup1992}. How can one estimate how much the orientation condition (ii) decreases the binding rate compared to $\ks$?

In the literature \cite{korennykh2007, zhou2010}, the orientation condition (ii) is sometimes accounted for by merely multiplying the rate constant $\ks$ by the product of the geometric correction factors,
\begin{align}\label{fsa}
\begin{split}
f_{A}:=\text{fraction of the $A$ protein surface area covered by binding sites},\\
f_{B}:=\text{fraction of the $B$ protein surface area covered by binding sites}.
\end{split}
\end{align}
The idea is that each protein collision has probability $f_{A}f_{B}\in(0,1)$ of having the binding sites aligned, and so the binding rate should be
\begin{align}\label{kgeo}
\kg
:=f_{A}f_{B}\ks.
\end{align}
However, experimentally measured protein-protein binding rates are typically a few orders of magnitude greater than the simple geometric estimate $\kg$ \cite{northrup1992} (note that binding sites usually occupy only a small portion of the protein surface, $f_{A}\ll1,f_{B}\ll1$ \cite{zhou2010}). While the estimate $\kg$ is simple and intuitive, it vastly underestimates the binding rate since, due to fine scale properties of Brownian motion, any proteins that collide once will collide many times in different orientations before they can diffuse away.

In this paper, we formulate and analyze a mathematical model of protein-protein binding to derive a bimolecular binding rate constant that includes both the proximity condition and the orientation condition given above. The model tracks a pair of diffusing spherical molecules, each with an arbitrary number of small binding sites on their surface, and the two molecules bind if and only if their binding sites come into contact (such molecules are often called ``patchy particles'' in the biochemistry literature \cite{zhang2004, fantoni2007, pawar2010, slyk2016, newton2015, roberts2014, klein2014}). Our analysis yields a first-principles derivation and estimate of both (a) how the orientation condition decreases the binding rate compared to $\ks$ in \eqref{ksmol} and (b) how purely diffusive processes increase the binding rate compared to $\kg$ in \eqref{kgeo}.

Mathematically, our model generalizes the {\black 1917} Smoluchowski model {\black \cite{smoluchowski1917}} in \eqref{pdeG}-\eqref{ksmol} to include the orientation condition (ii) given above. {\black In fact, our model generalizes the 1977 model of Berg and Purcell \cite{berg1977}, as our model reduces to their classical model if we assume one of the molecules is completely covered by binding sites.} By including {\black the orientation condition (ii)}, the equation \eqref{pdeG} becomes a high-dimensional ($\ge$ 7-dimensions), anisotropic diffusion equation, and the boundary condition \eqref{acG} at the reaction radius becomes a complicated mixed boundary condition. We apply formal matched asymptotic analysis \cite{lindsay2017} to this model in the case that the binding sites on each protein are small and well-separated (corresponding to a small surface area covered by binding sites, $f_{A}\ll1,f_{B}\ll1$). Our analysis yields a binding rate constant, $k=k_{0}$, which is much less than the Smoluchowski rate \eqref{ksmol} and much greater than the geometric estimate \eqref{kgeo} across a wide range of parameter values (that is, $\kg\ll k_{0}\ll\ks$). Our binding rate formula involves a dimensionless factor $\chi>0$ which is determined by the electrostatic capacitance of a certain 4-dimensional region embedded in 5-dimensions. As we do not have an exact analytical formula for $\chi$, we modify a recent kinetic Monte Carlo method \cite{bernoff2018b} to rapidly compute $\chi$ numerically. We then combine the quasi chemical formalism of {\v{S}}olc and Stockmayer \cite{solc1973} with recent asymptotic results \cite{lawley2019bp} to obtain a simple analytical approximation to $\chi$ which we show to be fairly accurate. This analysis further yields a binding rate formula that includes the effects of binding site competition/saturation. We verify our results by numerical simulations of the full system.

The rest of the paper is organized as follows. In section~\ref{results}, we describe the model and summarize our main results. In section~\ref{math}, we formulate the model more precisely and analyze the corresponding partial differential equation (PDE). In section~\ref{calculatingc0}, we develop a kinetic Monte Carlo method for computing $\chi$. In section~\ref{comp}, we apply the quasi chemical approximation. In section~\ref{numerical}, we verify our results by numerical simulations. We conclude by discussing applications and related work.

\section{Summary of main results}\label{results}

\begin{figure}
  \centering
    \includegraphics[width=0.5\textwidth]{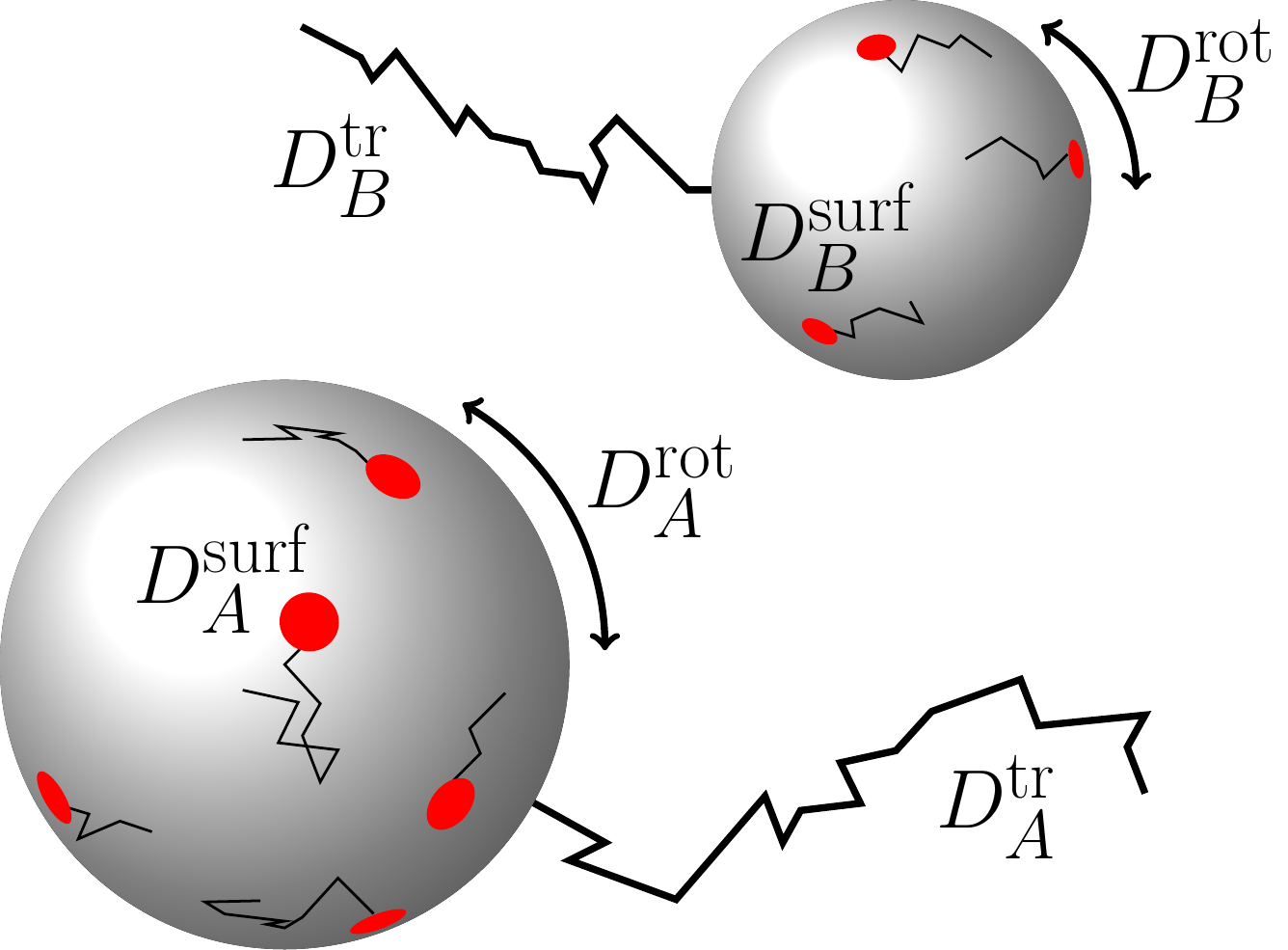}
 \caption{Two spherical molecules diffuse with respective translational diffusivities $\Dtr_{A}>0$ and $\Dtr_{B}>0$ and rotational diffusivities $\Drot_{A}\ge0$ and $\Drot_{B}\ge0$. The molecules have respectively $N_{A}\ge1$ and $N_{B}\ge1$ small, locally circular binding sites on their surfaces. The binding sites diffuse independently on the molecular surfaces with respective diffusivities $\Dsurf_{A}\ge0$ and $\Dsurf_{B}\ge0$. When the molecules come into contact, they bind if and only if a binding site on the $A$ molecule touches a binding site on the $B$ molecule; otherwise they reflect. }
 \label{schematicFig}
\end{figure}

Consider two spherical molecules, $A$ and $B$, with respective radii
\begin{align*}
R_{A}>0,\quad R_{B}>0,
\end{align*}
\emph{translational} diffusivities
\begin{align*}
\Dtr_{A}>0,\quad\Dtr_{B}>0,
\end{align*}
and \emph{rotational} diffusivities
\begin{align*}
\Drot_{A}\ge0,\quad\Drot_{B}\ge0.
\end{align*}
Suppose further that the molecules respectively have
\begin{align*}
N_{A}\ge1,\quad N_{B}\ge1
\end{align*}
small, locally circular binding sites on their surfaces with respective radii
\begin{align}\label{radii}
\eps a_{A}R_{A},\quad \eps a_{B}R_{B},
\end{align}
where $\eps\ll1$ is a small dimensionless parameter. The parameters $a_{A}$ and $a_{B}$ are order one dimensionless constants which allow the $A$ and $B$ binding sites to differ in size. We make no assumptions about the arrangements of the binding sites, except that they are well-separated, which means that the radii of the $A$ binding sites (respectively, $B$ binding sites) are much less than the typical distance between $A$ binding sites (respectively, $B$ binding sites). We further allow the possibility that the binding sites diffuse independently on the surfaces of their respective molecules with respective surface diffusivities
\begin{align*}
\Dsurf_{A}\ge0,\quad\Dsurf_{B}\ge0.
\end{align*}
To avoid trivial cases, we assume that the following ``effective'' diffusivities of the binding sites are strictly positive,
\begin{align*}
\Deff_{A}
&:=\Drot_{A}+R^{-2}\Dsurf_{A}>0,\quad
\Deff_{B}
:=\Drot_{B}+R^{-2}\Dsurf_{B}>0.
\end{align*}

Suppose the molecules bind if and only if a binding site on $A$ touches a binding site on $B$, otherwise they reflect. That is, the molecules bind if and only if they touch (proximity condition (i) above) \emph{and} the point of contact is in a binding site for \emph{both} molecules (orientation condition (ii) above). Note that if the binding requirement was merely that the point of contact is in a binding site for the $A$ molecule (and $\Drot_{A}=\Dsurf_{A}=0$), then we would obtain the 1977 model of Berg and Purcell \cite{berg1977}. Note also that if each molecule has only a single binding site ($N_{A}=N_{B}=1$), then we obtain the 1971 model of {\v{S}}olc and Stockmayer \cite{solc1971}. See Figure \ref{schematicFig} for a schematic representation of the model.

The initial state of the system can be described by a vector,
\begin{align}\label{is}
\begin{split}
&\Big((r,\theta_{0},\varphi_{0}),(\theta_{A}^{1},\varphi_{A}^{1}),\dots,(\theta_{A}^{N_{A}},\varphi_{A}^{N_{A}}),(\theta_{B}^{1},\varphi_{B}^{1}),\dots,(\theta_{B}^{N_{B}},\varphi_{B}^{N_{B}})\Big)\\
&\qquad\in[R_{A}+R_{B},\infty)\times\big([0,\pi]\times[0,2\pi)\big)^{1+N_{A}+N_{B}}\subset\R^{3+2N_{\f}+2N_{\m}},
\end{split}
\end{align}
which records the 3D location, $(r,\theta_{0},\varphi_{0})$, of the $B$ molecule relative to the $A$ molecule, as well as the 2D locations, $(\theta_{A}^{i},\varphi_{A}^{i})$ for $i\in\{1,\dots,N_{A}\}$ and $(\theta_{B}^{j},\varphi_{B}^{j})$ for $j\in\{1,\dots,N_{B}\}$, of the $N_{\f}+N_{\m}$ binding sites on the $A$ and $B$ molecules. Letting $p$ denote the probability that the two molecules never bind given the initial state \eqref{is}, we define the bimolecular binding rate constant $k_{0}>0$ analogously to \eqref{ksmol},
\begin{align}\label{k0def}
k_{0}
:=\Dtr\int_{r=R}\partial_{r}p\,\dd \Sigma,
\end{align}
where $\Dtr=\Dtr_{A}+\Dtr_{B}$ and the integration is over all possible initial states of the system \eqref{is} with $r=R:=R_{A}+R_{B}$ fixed at the reaction radius.

Using formal matched asymptotic analysis \cite{lindsay2017}, we show that (section~\ref{math})
\begin{align}\label{k0main}
k_{0}
\sim\eps^{3}N_{\f}N_{\m}\chi \ks,\quad\text{as }\eps\to0+,
\end{align}
where $\ks$ is the Smoluchowski rate \eqref{ksmol} and
\begin{align*}
\chi=\chi(\lambda_{A},\lambda_{B},a_{A},a_{B})>0
\end{align*}
is a dimensionless factor which depends on $a_{A}$, $a_{B}$, and the parameters 
\begin{align*}
\lambda_{A}
&:=\sqrt{1+\frac{R^{2}\Deff_\f}{\Dtr}}>1,
\quad
\lambda_{B}
:=\sqrt{1+\frac{R^{2}\Deff_\m}{\Dtr}}>1.
\end{align*}
That is, $\chi$ describes how the effective orientational diffusivities of $A$ and $B$ contribute to the binding rate. We modify a recent kinetic Monte Carlo method \cite{bernoff2018b} to rapidly and accurately compute $\chi$ (section~\ref{calculatingc0}).

In section~\ref{comp}, we combine recent asymptotic results \cite{lawley2019bp} for the case that one of the molecules is completely covered in binding sites with a heuristic quasi chemical approximation \cite{solc1973} to obtain a simple analytical approximation to $\chi$,
\begin{align}\label{chia}
\chiqc(\lambda_{A},\lambda_{B},a_{A},a_{B})
&:=\frac{a_{A}a_{B}(a_{A}\lambda_{B}+a_{B}\lambda_{A})}{4\pi}
\approx\chi(\lambda_{A},\lambda_{B},a_{A},a_{B}).
\end{align}
Using the kinetic Monte Carlo method of section~\ref{calculatingc0} we find that the relative error in the approximation \eqref{chia} is less than $16\%$ for $(R^{2}\Deff_{A}/\Dtr,R^{2}\Deff_{B}/\Dtr)\in[10^{-2},10]^{2}$ and $a_{A}=a_{B}=1$. Further, the analysis in section~\ref{comp} yields the following bimolecular binding rate formula which includes the effects of binding site competition/saturation,
\begin{align}\label{kintsmain}
\overline{k_{0}}:=\frac{\eps^{3}N_{A}N_{B}}{1/\chi+\eps^{2}\pi(N_{A}/(a_{B}\lambda_{B})+N_{B}/(\lambda_{A}a_{A}))+\eps^{3}N_{A}N_{B}}
\textcolor{black}{\ks\in(0,\ks)}.
\end{align}
In particular, \eqref{kintsmain} agrees with \eqref{k0main} in the limit $\eps\to0$ and has the correct limiting behavior if $N_{A}\to\infty$ and/or $\lambda_{A}\to\infty$ and/or $N_{B}\to\infty$ and/or $\lambda_{B}\to\infty$. That is,
\begin{align}\label{clb}
\begin{split}
\lim_{N_{B}\to\infty}\overline{k_{0}}
&=\lim_{\lambda_{B}\to\infty}\overline{k_{0}}
=\frac{\eps a_{A}\lambda_{A}N_{A}}{\pi+\eps a_{A}\lambda_{A}N_{A}}\ks
=:\overline{k_{A}},
\end{split}
\end{align}
where $\overline{k_{A}}$ is the binding rate derived in \cite{lawley2019bp} for the case that the $B$ molecule is completely covered in binding sites (of course, the analogous statement to \eqref{clb} holds if $N_{A}\to\infty$ and/or $\lambda_{A}\to\infty$). Note that since $\chi\approx\chiqc$, for simplicity we could replace $\chi$ by $\chiqc$ in the definition of $\overline{k_{0}}$ and obtain similar results. In section~\ref{numerical}, we compare our theoretical results to numerical simulations in order to (i) verify \eqref{k0main} and (ii) show that \eqref{kintsmain} is a good approximation to the bimolecular binding rate even away from the limits in \eqref{clb}.

\section{Matched asymptotic analysis} \label{math}

In this section, we derive the binding rate formula \eqref{k0main}. We begin by describing the stochastic binding model. 

\subsection{Stochastic problem formulation}\label{stochastic}

Consider first the case of zero rotational diffusion. That is, suppose that $\Drot_{A}=\Drot_{B}=0$ and $\Dsurf_{A}>0$, $\Dsurf_{B}>0$. We show below that our results are quickly extended to the general case $\Drot_{A}\ge0$, $\Drot_{B}\ge0$, $\Dsurf_{A}\ge0$, $\Dsurf_{B}\ge0$ with $\Deff_{A}:=\Drot_{A}+R^{-2}\Dsurf_{A}>0$, $\Deff_{B}:=\Drot_{B}+R^{-2}\Dsurf_{B}>0$.

 Fixing the reference frame on the $A$ molecule, the state of the system at time $t\ge0$ can be described by the 3D position in spherical coordinates of the center of the spherical $B$ molecule,
\begin{align}\label{0c}
(X(t),\Theta_{0}(t),\Phi_{0}(t))\in[R,\infty)\times[0,\pi]\times[0,2\pi),
\end{align}
and the 2D positions of the $A$ and $B$ binding sites. We denote the spherical angular coordinates of the center of the $i$th $A$ binding site at time $t\ge0$ by
\begin{align}\label{Ac}
(\Theta_{\f}^{i}(t),\Phi_{\f}^{i}(t))\in[0,\pi]\times[0,2\pi),\quad i\in\{1,\dots,N_{A}\}.
\end{align}
Rather than tracking the centers of the $B$ binding sites, it is convenient to track the positions of their antipodal points on the $B$ molecule at time $t\ge0$, which we denote by 
\begin{align}\label{Bc}
(\Theta_{\m}^{i}(t),\Phi_{\m}^{i}(t))\in[0,\pi]\times[0,2\pi),\quad i\in\{1,\dots,N_{B}\}.
\end{align}
Naturally, the coordinates in \eqref{0c} and \eqref{Ac} take the center of the $A$ molecule to be their origin, while the coordinates in \eqref{Bc} take the center of the $B$ molecule as their origin. All three sets of coordinates use the same $z$-direction to define their north poles.

Since the molecules have respective translational diffusivities $\Dtr_{A}$ and $\Dtr_{B}$, it follows that the coordinates in \eqref{0c} satisfy the following stochastic differential equations (SDEs) with $\Dtr:=\Dtr_{A}+\Dtr_{B}$,
\begin{align}\label{pSDE}
\begin{split}
\dd X(t)
&=\frac{2\Dtr}{X(t)}\dd t+\sqrt{2\Dtr}\,\dd W_{X}(t),\\
\dd\Theta_{0}(t)
&=\frac{\Dtr}{(X(t))^{2}\tan(\Theta_{0}(t))}\,\dd t
+\frac{\sqrt{2\Dtr}}{X(t)}\,\dd W_{\Theta_{0}}(t),\\
\dd\Phi_{0}(t)
&=\frac{\sqrt{2\Dtr}}{X(t)\sin(\Theta_{0}(t))}\,\dd W_{\Phi_{0}}(t),
\end{split}
\end{align}
where $W_{X}$, $W_{\Theta_{0}}$, and $W_{\Phi_{0}}$ are independent standard Brownian motions. The coordinates in \eqref{Ac} satisfy
\begin{align}\label{ASDE}
\begin{split}
\dd\Theta_{A}^{i}(t)
&=\frac{\Deff_{A}}{\tan(\Theta_{A}^{i}(t))}\,\dd t
+\sqrt{2\Deff_{A}}\,\dd W_{\Theta_{A}^{i}}(t),\\
\dd\Phi_{A}^{i}(t)
&=\frac{\sqrt{2\Deff_{A}}}{\sin(\Theta_{A}^{i}(t))}\,\dd W_{\Phi_{A}^{i}}(t),\quad i\in\{1,\dots,N_{A}\},
\end{split}
\end{align}
where $W_{\Theta_{A}^{i}}$ and $W_{\Phi_{A}^{i}}$ are independent standard Brownian motions. Similarly, the coordinates in \eqref{Bc} satisfy
\begin{align}\label{BSDE}
\begin{split}
\dd\Theta_{B}^{i}(t)
&=\frac{\Deff_{B}}{\tan(\Theta_{B}^{i}(t))}\,\dd t
+\sqrt{2\Deff_{B}}\,\dd W_{\Theta_{B}^{i}}(t),\\
\dd\Phi_{B}^{i}(t)
&=\frac{\sqrt{2\Deff_{B}}}{\sin(\Theta_{B}^{i}(t))}\,\dd W_{\Phi_{B}^{i}}(t),\quad i\in\{1,\dots,N_{B}\},
\end{split}
\end{align}
where $W_{\Theta_{B}^{i}}$ and $W_{\Phi_{B}^{i}}$ are independent standard Brownian motions. Note that $\Dtr$ has units of length squared per time, whereas $\Deff_{A}$ and $\Deff_{B}$ have units of inverse time.

\begin{figure}
  \centering
    \includegraphics[width=1\textwidth]{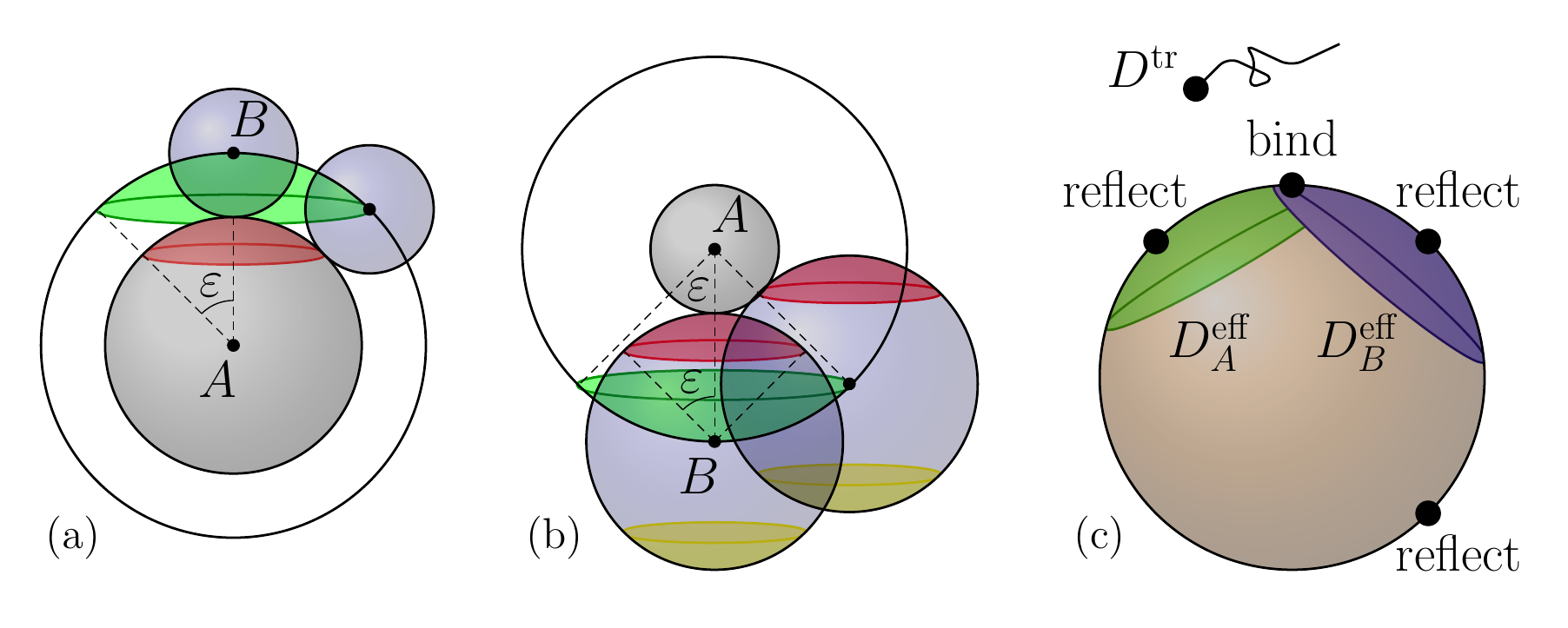}
 \caption{(a) The $B$ molecule (blue sphere) with radius $R_{B}$ touches the $A$ molecule (grey sphere) with radius $R_{A}$ if and only if their centers are distance $R=R_{A}+R_{B}$ apart. An $A$ binding site (red region) is a spherical cap with polar angle $\eps$. The $B$ molecule hits an $A$ binding site if and only if the center of the $B$ molecule hits a spherical cap (green region) with angle $\eps$ on the sphere of radius $R$. (b) The red region is a $B$ binding site and the yellow region is its antipodal region on the $B$ molecule. The green region is the same cap as the yellow region, but placed on the sphere of radius $R$. The $A$ molecule (grey sphere) hits a $B$ binding site (red region) if and only if the center of the $B$ molecule (blue sphere) is in the green region on the sphere of radius $R$. (c) The $A$ and $B$ binding sites are in contact if and only if the center of the $B$ molecule (black dot) is on the sphere of radius $R$ in the intersection of the spherical caps of the $A$ binding site (green region) and the antipodal $B$ binding site (blue region). The black dot diffuses with translational diffusivity $\Dtr$ and the green and blue regions diffuse independently on the surface of the sphere with respective diffusivities $\Deff_{A}$ and $\Deff_{B}$. In this figure, {\black each sphere has one binding site, $N_{A}=N_{B}=1$}.}
 \label{figschem2}
\end{figure}

For a pair of angular spherical coordinates, $(\theta',\varphi')\in[0,\pi]\times[0,2\pi)$, and a polar angle $\eps\in(0,\pi/2]$, define the spherical cap
\begin{align*}
\Gamma(\theta',\varphi',\eps)
:=\big\{(\theta,\varphi)\in[0,\pi]\times[0,2\pi):(\theta-\theta')^{2}+\sin^{2}(\theta')(\varphi-\varphi')^{2}<\eps^{2}\big\}.
\end{align*}
Each $A$ binding site is the spherical cap on the $A$ molecule centered at \eqref{Ac} with polar angle $\eps a_{\f}$. Similarly, each $B$ binding site is the spherical cap on the $B$ molecule centered at the point that is antipodal to \eqref{Bc} with polar angle $\eps a_{\m}$. Figure~\ref{figschem2} illustrates this geometry for the case {\black of a single binding site on each sphere, $N_{A}=N_{B}=1$}.

It is readily apparent (see Figure~\ref{figschem2}) that the $i$th $A$ binding site and the $j$th $B$ binding site are in contact at time $t\ge0$ if and only if $X(t)=R$ and the angular position of the $B$ molecule is in the intersection of the two spherical caps,
\begin{align*}
(\Theta_{0}(t),\Phi_{0}(t))\in \Gamma(\Theta_{\f}^{i}(t),\Phi_{\f}^{i}(t),\eps a_{\f})\cap\Gamma(\Theta_{\m}^{j}(t),\Phi_{\m}^{j}(t),\eps a_{\m}).
\end{align*}
It follows that this problem is equivalent to (i) a set of $N_{A}$ spherical caps and a set of $N_{B}$ spherical caps which all diffuse independently on the surface of a sphere with radius $R$ and (ii) a point particle at position $(X(t),\Theta_{0}(t),\varphi_{0}(t))$ that diffuses exterior to this sphere and is absorbed at the sphere if and only if it hits the intersection of an $A$ spherical cap with a $B$ spherical cap (otherwise it reflects from the sphere). In particular, the $A$ and $B$ molecules bind if and only if this point particle reaches the intersection of these two sets of spherical caps (see Figure~\ref{figschem2}(c)).

Let $\tau\ge0$ denote the random time when the two molecules bind,
\begin{align}\label{tau}
\tau
&:=\inf\Big\{t>0:X(t)=R,\,(\Theta_{0}(t),\Phi_{0}(t))\in\Lambda_{\f}(\Theta_{A}(t),\Phi_{A}(t))\cap\Lambda_{\m}(\Theta_{B}(t),\Phi_{B}(t))\Big\},
\end{align}
where we have defined the unions of the $A$ and $B$ caps respectively as
\begin{align*}
\Lambda_\f(\Theta_{A}(t),\Phi_{A}(t))
&:= \cup_{i=1}^{N_\f} \Gamma(\Theta_{\f}^{i}(t),\Phi_{\f}^{i}(t),\eps a_{\f}),
\\
\Lambda_\m(\Theta_{\m}(t),\Phi_{\m}(t))
&:= \cup_{i=1}^{N_\m} \Gamma(\Theta_{\m}^{i}(t),\Phi_{\m}^{i}(t),\eps a_{\m}),
\end{align*}
and the vectors of angles,
\begin{align*}
\Theta_{A}(t)
&:=(\Theta_{A}^{1}(t),\dots,\Theta_{A}^{N_{A}}(t))\in[0,\pi]^{N_{A}},\\
\Phi_{A}(t)
&:=(\Phi_{A}^{1}(t),\dots,\Phi_{A}^{N_{A}}(t))\in[0,2\pi)^{N_{A}},\\
\Theta_{B}(t)
&:=(\Theta_{B}^{1}(t),\dots,\Theta_{B}^{N_{B}}(t))\in[0,\pi]^{N_{B}},\\
\Phi_{B}(t)
&:=(\Phi_{B}^{1}(t),\dots,\Phi_{B}^{N_{B}}(t))\in[0,2\pi)^{N_{B}}.
\end{align*}
Let $p(r,\theta_{0},\varphi_{0},\theta_{\f},\varphi_{\f},\theta_{\m},\varphi_{\m})$ denote the probability that the molecules never bind, conditioned on the initial state of the system,
\begin{align}
\begin{split}\label{pdef}
p(r,\theta_{0},\varphi_{0},\theta_{\f},\varphi_{\f},\theta_{\m},\varphi_{\m})
&:=\P\big(\tau=\infty\,\big|\,X(0)=r,\Theta_{0}(0)=\theta_{0},\Phi_{0}(0)=\varphi_{0},\\
&\qquad\Theta_{\f}(0)=\theta_{\f},\Phi_{\f}(0)=\varphi_{\f},\Theta_{\m}(0)=\theta_{\m},\Phi_{\m}(0)=\varphi_{\m}\big),
\end{split}
\end{align}
where the arguments of the function are the initial state of the system, $r\in[R,\infty)$,
\begin{align*}
(\theta_{0},\varphi_{0})
&\in[0,\pi]\times[0,2\pi),\\
(\theta_{A},\varphi_{A})
=((\theta_{A}^{1},\varphi_{A}^{1}),\dots,(\theta_{A}^{N_{A}},\varphi_{A}^{N_{A}}))
&\in([0,\pi]\times[0,2\pi))^{N_{A}},\\
(\theta_{B},\varphi_{B})
=((\theta_{B}^{1},\varphi_{B}^{1}),\dots,(\theta_{B}^{N_{B}},\varphi_{B}^{N_{B}}))
&\in([0,\pi]\times[0,2\pi))^{N_{B}}.
\end{align*}

\subsection{PDE boundary value problem}

Define the elliptic operators
\begin{align}\label{la}
\LL_{A}:=\sum_{i=1}^{N_{A}}\L_{A}^{i},\quad
\LL_{B}:=\sum_{i=1}^{N_{B}}\L_{B}^{i},
\end{align}
where $\L_\f^{i}$ denotes the Laplace-Beltrami operator acting on $(\theta_\f^{i},\varphi_\f^{i})$,
\begin{align}\label{lb}
\L_{A}^{i}:=(\sin(\theta_{A}^{i}))^{-2}\partial_{\varphi_{A}^{i}\varphi_{A}^{i}}+\cot(\theta_{A}^{i})\partial_{\theta_{A}^{i}}+\partial_{\theta_{A}^{i}\theta_{A}^{i}},\quad i\in\{1,\dots,N_{A}\},
\end{align}
and similarly for $\L_{B}^{i}$. Let $\Delta_{0}$ denote the Laplacian acting on $(r,\theta_{0},\varphi_{0})$,
\begin{align}\label{l0}
\Delta_{0}=\frac{2}{r}\partial_{r}+\partial_{rr}+\frac{1}{r^{2}}\L_{0},
\end{align}
where $\L_{0}$ acts on $(\theta_{0},\varphi_{0})$ as in \eqref{lb}.

It is straightforward to show that $p$ satisfies the following elliptic PDE,
\begin{align}\label{ppde}
0=(\Dtr\Delta_{0}+\Deff_\f \LL_\f+\Deff_\m \LL_\m)p,\quad r>R.
\end{align}
Since molecules starting far from each other will never bind, we obtain the following far-field condition which is identical to \eqref{ffc},
\begin{align}\label{pff}
\lim_{r\to\infty}p=1.
\end{align}
Finally, since the molecules bind if the binding sites are in contact and otherwise reflect, we obtain the following mixed boundary conditions at the reaction radius $r=R$,
\begin{align}\label{pbc}
\begin{split}
p
&=0,\quad r=R,\,(\theta_{0},\varphi_{0})\in \Lambda(\theta_{\f},\varphi_{\f})\cap\Lambda(\theta_{\m},\varphi_{\m}),\\
\partial_{r}p
&=0,\quad r=R,\,(\theta_{0},\varphi_{0})\notin \Lambda(\theta_{\f},\varphi_{\f})\cap\Lambda(\theta_{\m},\varphi_{\m}).
\end{split}
\end{align}
The PDE boundary value problem in \eqref{ppde}-\eqref{pbc} generalizes the classical Smoluchowski model in \eqref{pdeG}-\eqref{acG}.

\subsection{Outer expansion}

We now apply formal matched asymptotic analysis to the $(3+2N_{A}+2N_{B})$-dimensional PDE boundary value problem in \eqref{ppde}-\eqref{pbc}. Our approach follows the methods employed in \cite{lindsay2017} to analyze a similar 3D problem. These formal methods are related to the strong localized perturbation analysis pioneered in \cite{ward93,ward93b}.

It is straightforward to show that $p$ has the following behavior at far-field,
\begin{align}\label{pffc}
1-p\sim \frac{C}{r},\quad r\to\infty,
\end{align}
for some constant $C\in(0,R)$. In an analogy to electrostatics, we refer to $C$ as the capacitance. It follows that the bimolecular binding rate in \eqref{k0def} is related to the capacitance by
\begin{align}\label{k0C}
k_{0}=\frac{C}{R}\ks.
\end{align}

It is convenient to work with the following rescaling of $p$,
\begin{align}\label{vdef}
v:=\frac{-p}{C}.
\end{align}
We expect that (i) $v$ has boundary layers near the absorbing boundary conditions and (ii) $C=\O(\eps^{3})$ as $\eps\to0$. We thus introduce the outer expansion,
\begin{align}\label{outer}
v\sim
\eps^{-3}v_{0}+v_{1}+\cdots,
\end{align}
where $v_{0}$ is a constant and $v_{1}$ is a function. The outer expansion \eqref{outer} is valid away from the boundary layers. {\black We note that one can derive the scaling $C=\O(\eps^{3})$ by considering the case of zero rotational and surface diffusion analyzed in section~\ref{zero} below.}

Using the definition of $v$ in \eqref{vdef} and plugging the outer expansion \eqref{outer} into \eqref{ppde} implies that $v_{1}$ satisfies
\begin{align}\label{v1first}
\begin{split}
0
&=(\Dtr\Delta_{0}+\Deff_\f \LL_\f+\Deff_\m \LL_\m)v_{1},\quad r>R,\\
\partial_{r}v_{1}
&=0,\quad r=R,\,(\theta_{0},\varphi_{0})\notin \big\{\cup_{i=1}^{N_{A}}\{(\theta_\f^{i},\varphi_\f^{i})\}\big\} \cap \big\{\cup_{j=1}^{N_{B}}\{(\theta_\m^{j},\varphi_\m^{j})\}\big\} .
\end{split}
\end{align}
Notice that the binding sites have shrunk to points from the perspective of the outer solution $v_{1}$.

\subsection{Inner expansion}

We now determine the singular behavior of $v_{1}$ as
\begin{align*}
(r,\theta_\f^{i} ,\varphi_\f^{i} ,\theta_\m^{j} ,\varphi_\m^{j} )\to(R,\theta_{0},\varphi_{0},\theta_{0},\varphi_{0}),\quad\text{for some }i\in\{1,\dots,N_{A}\},\,j\in\{1,\dots,N_{B}\}.
\end{align*}
First, introduce the stretched coordinates,
\begin{align}\label{stretch}
\begin{split}
{{z}}
&:=\eps^{-1}(\tfrac{r}{R}-1),\\
t_\f
&:=\eps^{-1}(\theta_\f^{i}-\theta_{0}),\\
p_\f
&:=\eps^{-1}\sin(\theta_{0})(\varphi_\f^{i}-\varphi_{0}),\\
t_\m 
&:=\eps^{-1}(\theta_\m^{j} -\theta_{0}),\\
p_\m 
&:=\eps^{-1}\sin(\theta_{0})(\varphi_\m^{j} -\varphi_{0}).
\end{split}
\end{align}
We note that $t_{A},p_{A}$ depend on $i$ and $t_{B},p_{B}$ depend on $j$, but we suppress this dependence to simplify notation. Next, introduce the linear combinations,
\begin{alignat*}{2}
{{x_{\f}}}
&=c_{11}t_\f +c_{12}t_\m,\\
{{y_{\f}}}
&=d_{11}p_\f +d_{12}p_\m,\\
{{x_{\m}}}
&=c_{21}t_\f +c_{22}t_\m,\\
{{y_{\m}}}
&=d_{21}p_\f +d_{22}p_\m,
\end{alignat*}
where the 8 constants,
\begin{align}\label{constants}
c_{11},c_{12},c_{21},c_{22},d_{11},d_{12},d_{21},d_{22}
\end{align}
are yet to be determined. We then define the inner solution
\begin{align}\label{innersol0}
\begin{split}
w({{z}},{{x_{\f}}},{{y_{\f}}},{{x_{\m}}},{{y_{\m}}})
&=w\Big({{z}},{{x_{\f}}},{{y_{\f}}},{{x_{\m}}},{{y_{\m}}};(\theta_{0},\varphi_{0}),\\
&\qquad(\theta_{A}^{1},\varphi_{A}^{1}),\dots,(\theta_{A}^{i-1},\varphi_{A}^{i-1}),(\theta_{A}^{i+1},\varphi_{A}^{i+1}),\dots,(\theta_{A}^{N_{A}},\varphi_{A}^{N_{A}}))\\
&\qquad(\theta_{B}^{1},\varphi_{B}^{1}),\dots,(\theta_{B}^{j-1},\varphi_{B}^{j-1}),(\theta_{B}^{j+1},\varphi_{B}^{j+1}),\dots,(\theta_{B}^{N_{B}},\varphi_{B}^{N_{B}})\Big)\\
&:=v\Big(R+\eps R {{z}},\theta_{0},\varphi_{0},
(\theta_{A}^{1},\varphi_{A}^{1}),\dots,(\theta_{A}^{i-1},\varphi_{A}^{i-1}),\\
&\qquad\theta_{0}+\eps\Big(\frac{c_{12}x_{B}-c_{22}x_{A}}{c_{12}c_{21}-c_{11}c_{22}}\Big),
\varphi_{0}+\eps\Big(\frac{d_{12}y_{B}-d_{22}y_{A}}{\sin(\theta_{0})(d_{12}d_{21}-d_{11}d_{22})}\Big),\\
&\qquad(\theta_{A}^{i+1},\varphi_{A}^{i+1}),\dots,(\theta_{A}^{N_{A}},\varphi_{A}^{N_{A}}),
(\theta_{B}^{1},\varphi_{B}^{1}),\dots,(\theta_{B}^{j-1},\varphi_{B}^{j-1}),\\
&\qquad\theta_{0}+\eps\Big(\frac{c_{21}x_{A}-c_{11}x_{B}}{c_{12}c_{21}-c_{11}c_{22}}\Big),
\varphi_{0}+\eps\Big(\frac{d_{21}y_{A}-d_{11}y_{B}}{\sin(\theta_{0})(d_{12}d_{21}-d_{11}d_{22})}\Big)\\
&\qquad(\theta_{B}^{j+1},\varphi_{B}^{j+1}),\dots,(\theta_{B}^{N_{B}},\varphi_{B}^{N_{B}})\Big).
\end{split}
\end{align}
In words, the inner solution $w$ zooms in on $(r,\theta_\f^{i} ,\varphi_\f^{i} ,\theta_\m^{j} ,\varphi_\m^{j} )$ near $(R,\theta_{0},\varphi_{0},\theta_{0},\varphi_{0})$.

We now choose the constants in \eqref{constants} so that the inner solution satisfies an isotropic diffusion equation to leading order as $\eps\to0$. The definition of the inner solution in \eqref{innersol0} implies that
\begin{align*}
&v(r,\theta_{0},\varphi_{0},\theta_\f,\varphi_\f,\theta_\m ,\varphi_\m)\\
&=w(
\eps^{-1}(r/R-1),
c_{11}t_\f+c_{12}t_\m ,     d_{11}p_\f+d_{12}p_\m ,\\
&\quad\qquad\qquad\qquad\qquad c_{21}t_\f+c_{22}t_\m ,     d_{21}p_\f+d_{22}p_\m),\\
&=w(
\eps^{-1}(r/R-1),
c_{11}\eps^{-1}(\theta_\f^{i}-\theta_{0})+c_{12}\eps^{-1}(\theta_\m^{j} -\theta_{0}),\\
&\quad\qquad\qquad\qquad\qquad     d_{11}\eps^{-1}\sin(\theta_{0})(\varphi_\f^{i}-\varphi_{0})+d_{12}\eps^{-1}\sin(\theta_{0})(\varphi_\m^{j} -\varphi_{0}),\\
&\quad\qquad\qquad\qquad\qquad c_{21}\eps^{-1}(\theta_\f^{i}-\theta_{0})+c_{22}\eps^{-1}(\theta_\m^{j} -\theta_{0}),\\
&\quad\qquad\qquad\qquad\qquad     d_{21}\eps^{-1}\sin(\theta_{0})(\varphi_\f^{i}-\varphi_{0})+d_{22}\eps^{-1}\sin(\theta_{0})(\varphi_\m^{j} -\varphi_{0})).
\end{align*}
We now calculate the leading order terms in the differential operator,
\begin{align*}
\mathbb{L}:={\Dtr}\Delta_{0}+\Deff_\f\LL_\f+\Deff_\m \LL_\m, 
\end{align*}
where $\Delta_{0}$, $\LL_{A}$, and $\LL_{B}$ are defined in \eqref{la} and \eqref{l0}.
Specifically, 
\begin{align*}
\eps^{2}\mathbb{L}v
&=
\O(\eps)+R^{-2}{\Dtr}w_{{{z}}{{z}}}
+R^{-2}{\Dtr}\Big\{(-c_{11}-c_{12})^{2}w_{{{x_{\f}}}{{x_{\f}}}}
+(-c_{21}-c_{22})^{2}w_{{{x_{\m}}}{{x_{\m}}}}\\
&\quad
+2(-c_{11}-c_{12})(-c_{21}-c_{22})w_{{{x_{\f}}}{{x_{\m}}}}\\
&\quad
+(-d_{11}-d_{12})^{2}w_{{{y_{\f}}}{{y_{\f}}}}
+(-d_{21}-d_{22})^{2}w_{{{y_{\m}}}{{y_{\m}}}}\\
&\quad
+2(-d_{11}-d_{12})(-d_{21}-d_{22})w_{{{y_{\f}}}{{y_{\m}}}}
\Big\}\\
&\quad
+\Deff_\f\Big\{c_{11}^{2}w_{{{x_{\f}}}{{x_{\f}}}}+c_{21}^{2}w_{{{x_{\m}}}{{x_{\m}}}}
+2c_{11}c_{21}w_{{{x_{\f}}}{{x_{\m}}}}\\
&\quad
+d_{11}^{2}w_{{{y_{\f}}}{{y_{\f}}}}+d_{21}^{2}w_{{{y_{\m}}}{{y_{\m}}}}
+2d_{11}d_{21}w_{{{y_{\f}}}{{y_{\m}}}}\Big\}\\
&\quad
+\Deff_\m \Big\{c_{12}^{2}w_{{{x_{\f}}}{{x_{\f}}}}+c_{22}^{2}w_{{{x_{\m}}}{{x_{\m}}}}
+2c_{12}c_{22}w_{{{x_{\f}}}{{x_{\m}}}}\\
&\quad
+d_{12}^{2}w_{{{y_{\f}}}{{y_{\f}}}}+d_{22}^{2}w_{{{y_{\m}}}{{y_{\m}}}}
+2d_{12}d_{22}w_{{{y_{\f}}}{{y_{\m}}}}\Big\},
\end{align*}
which upon collecting terms shows that $w$ satisfies the following leading order equation as $\eps\to0$,
\begin{align*}
0=
w_{zz}
&
+\Big\{(c_{11}+c_{12})^{2}+D_\f c_{11}^{2}+D_\m c_{12}^{2}\Big\}w_{{{x_{\f}}}{{x_{\f}}}}\\
&
+\Big\{(d_{11}+d_{12})^{2}+D_\f d_{11}^{2}+D_\m d_{12}^{2}\Big\}w_{{{y_{\f}}}{{y_{\f}}}}\\
&
+\Big\{(c_{21}+c_{22})^{2}+D_\f c_{21}^{2}+D_\m c_{22}^{2}\Big\}w_{{{x_{\m}}}{{x_{\m}}}}\\
&
+\Big\{(d_{21}+d_{22})^{2}+D_\f d_{21}^{2}+D_\m d_{22}^{2}\Big\}w_{{{y_{\m}}}{{y_{\m}}}}\\
&
+2\Big\{(c_{11}+c_{12})(c_{21}+c_{22})+D_\f c_{11}c_{21}+D_\m c_{12}c_{22}\Big\}w_{{{x_{\f}}}{{x_{\m}}}}\\
&
+2\Big\{(d_{11}+d_{12})(d_{21}+d_{22})+D_\f d_{11}d_{21}+D_\m d_{12}d_{22}\Big\}w_{{{y_{\f}}}{{y_{\m}}}}+\O(\eps),
\end{align*}
where we have defined the ratios
\begin{align}\label{Ddef}
D_{A}:=R^{2}\Deff_\f/\Dtr,
\quad
D_{B}:=R^{2}\Deff_\m/\Dtr.
\end{align}

We now choose the constants in \eqref{constants} so that all the pure second derivative terms have the same coefficient and the mixed partial derivative terms vanish. In particular, $c_{11},c_{12},c_{21},c_{22}$ must satisfy
\begin{align}\label{nu}
\begin{split}
1
&=(c_{11}+c_{12})^{2}+D_\f  c_{11}^{2}+D_\m c_{12}^{2}\\
1
&=(c_{21}+c_{22})^{2}+D_\f c_{21}^{2}+D_\m c_{22}^{2}\\
0
&=(c_{11}+c_{12})(c_{21}+c_{22})+D_\f c_{11}c_{21}+D_\m c_{12}c_{22},
\end{split}
\end{align}
and $d_{11},d_{12},d_{21},d_{22}$ must satisfy the same equations. We thus have 6 nonlinear equations for the 8 constants in \eqref{constants}. There is not a unique solution. Nevertheless, we choose the following solution,
\begin{align}
\label{cdefs}
\begin{split}
c_{11}
&:=d_{11}
:=
 \lambda_\f ^{-1},\quad\quad\quad
 c_{12}
:=d_{12}
:=0,\\
c_{21}
:=d_{21}
&:=
\lambda_\f ^{-1}\frac{1}{\sqrt{\lambda_\f ^{2}\lambda_\m ^{2}-1}},\quad
c_{22}
:=d_{22}
:=
-\lambda_\f \frac{1}{\sqrt{\lambda_\f ^{2}\lambda_\m ^{2}-1}},
\end{split}
\end{align}
where
\begin{align}\label{lambdadef}
\lambda_\f := \sqrt{1+D_{A}},
\quad
\lambda_\m  := \sqrt{1+D_{B}}.
\end{align}
{\black (We have not investigated all the solutions to \eqref{nu}, but we have investigated some other solutions and found that the asymptotic behavior of the binding rate $k_{0}$ in \eqref{k0def} that ultimately results does not depend on the choice of solution to \eqref{nu}.)} By the choice \eqref{cdefs}, we note that
\begin{align}
\label{angledefs0}
\begin{split}
\theta_\f 
&= \theta_{0}+\eps\frac{{{x_{\f}}}}{c_{11}},\\
\varphi_\f 
&=\varphi_{0}+\eps\frac{{{y_{\f}}} }{c_{11} \sin(\theta_{0})}, \\
\theta_\m 
&=\theta_{0}+
\eps\frac{  (c_{11} {{x_{\m}}}-c_{21} {{x_{\f}}})}{c_{11}c_{22}}, \\
\varphi_\m 
&=\varphi_{0}+\eps\frac{ (c_{11} {{y_{\m}}}-c_{21} {{y_{\f}}})}{c_{11} c_{22} \sin(\theta_{0})}.
\end{split}
\end{align}

By construction, the inner solution in \eqref{innersol0} is harmonic in the 5 variables, $({{z}},{{x_{\f}}},{{y_{\f}}},{{x_{\m}}},{{y_{\m}}})$, to leading order. Indeed, if we introduce the inner expansion,
\begin{align*}
w
\sim\eps^{-3}w_{0}+\cdots,
\end{align*}
then the calculation above implies that $w_{0}$ is harmonic in half of 5-dimensional space, $({{z}},{{x_{\f}}},{{y_{\f}}},{{x_{\m}}},{{y_{\m}}})\in(0,\infty)\times\R^{4}$,
\begin{align}\label{harmonic}
\big(\partial_{{{z}}{{z}}}+
\partial_{{{x_{\f}}}{{x_{\f}}}}+
\partial_{{{y_{\f}}}{{y_{\f}}}}+
\partial_{{{x_{\m}}}{{x_{\m}}}}+
\partial_{{{y_{\m}}}{{y_{\m}}}}\big)w_{0}=0,\quad {{z}}>0,\,({{x_{\f}}},{{y_{\f}}},{{x_{\m}}},{{y_{\m}}})\in\R^{4}.
\end{align}
Furthermore, the boundary conditions at ${{z}}=0$, are 
\begin{align}\label{w0bc}
\begin{split}
w_{0}
&=0,\quad {{z}}=0,\; {{x_{\f}}}^{2}+{{y_{\f}}}^{2}<c_{11}^{2} a_\f^2,\; \Big({{x_{\m}}}-\frac{c_{21}}{c_{11}}{{x_{\f}}}\Big)^{2}+\Big({{y_{\m}}}-\frac{c_{21}}{c_{11}}{{y_{\f}}}\Big)^{2}<c_{22}^{2} a_\m ^2\\
\partial_{{{z}}}w_{0}
&=0,\quad {{z}}=0,\; \text{otherwise}.
\end{split}
\end{align}
In words, $w_{0}=0$ if the following 3 conditions are satisfied simultaneously, (i) ${{z}}=0$, (ii) $({{x_{\f}}},{{y_{\f}}})$ is in a disk of radius $c_{11} a_\f >0$ centered at $(0,0)$, and (iii) $({{x_{\m}}},{{y_{\m}}})$ is in a disk of radius $c_{22} a_\m  >0$ centered at $\frac{c_{22}}{c_{11}}({{x_{\f}}},{{y_{\f}}})$. Otherwise, $\partial_{{{z}}}w_{0}=0$ if ${{z}}=0$.

To derive \eqref{w0bc}, note that if $r=R$ and $(\theta_0, \varphi_0) \in \Gamma (\theta_\f^{i}, \varphi_\f^{i}, \eps a_{A}) \cap \Gamma(\theta_\f^{j}, \varphi_\f^{j}, \eps a_{B})$, 
then $v=0$. By definition, $(\theta_0, \varphi_0) \in \Gamma (\theta_\f^{i}, \varphi_\f^{i}, \eps a_{A})$ means
\begin{align*}
(\theta_0 - \theta_\f^{i})^2 + \sin^2(\theta_\f^{i})(\varphi_0 - \varphi_\f^{i})^2 < \eps^2 a_\f^2,
\end{align*}
and using \eqref{angledefs0}, we find that this implies
\begin{align}\label{lowo}
\eps^2 \frac{x_{A}^2}{c_{11}^2} + \frac{\sin^2(\theta_0 + \O(\eps)) \eps^2 y_{A}^2}{ c_{11}^2 \sin^2(\theta_0)} < \eps^2 a_\f^2.
\end{align}
Taking terms of lowest order in $\eps$ in \eqref{lowo} and simplifying yields the condition ${{x_{\f}}}^{2}+{{y_{\f}}}^{2}<c_{11}^{2} a_\f^2$ in \eqref{w0bc}. The condition $({{x_{\m}}}-\frac{c_{21}}{c_{11}}{{x_{\f}}})^{2}+({{y_{\m}}}-\frac{c_{21}}{c_{11}}{{y_{\f}}})^{2}<c_{22}^{2} a_\m ^2$ in \eqref{w0bc} is obtained similarly.

\subsection{Matching}

It follows from electrostatics \cite{jackson1975} that $w_{0}$ has the far-field behavior
\begin{align}\label{w0far}
w_{0}\sim \alpha\Big(1-\frac{c_{0}(c_{11},c_{21},c_{22},a_\f,a_\m )}{\rho^{3}}\Big),\;\text{as }
\rho:=\sqrt{{{z}}^{2}+{{x_{\f}}}^{2}+{{y_{\f}}}^{2}+{{x_{\m}}}^{2}+{{y_{\m}}}^{2}}\to\infty,
\end{align}
where $\alpha$ is a constant to be determined by matching to the outer solution, and
\begin{align*}
c_{0}=c_{0}(D_{A},D_{B},a_{A},b_{B})
\end{align*}
is a constant depending on $D_{A},D_{B},a_{A},b_{B}$. In particular, $c_{0}$ is the electrostatic capacitance of the following 4D region embedded in $\R^{5}$,
\begin{align}\label{region}
\begin{split}
\mathcal{R}
&:=\Big\{(z,{{x_{\f}}},{{y_{\f}}},{{x_{\m}}},{{y_{\m}}})\in\R^{5}:z=0,\,{{x_{\f}}}^{2}+{{y_{\f}}}^{2}<c_{11}^{2} a_\f^2,\\
&\qquad\qquad\qquad\qquad\qquad\qquad\Big({{x_{\m}}}-\frac{c_{21}}{c_{11}}{{x_{\f}}}\Big)^{2}+\Big({{y_{\m}}}-\frac{c_{21}}{c_{11}}{{y_{\f}}}\Big)^{2}<c_{22}^{2} a_\m ^2\Big\}.
\end{split}
\end{align}
For the remainder of this section, we carry out our calculations in terms of $c_{0}$. In section~\ref{calculatingc0}, we develop a numerical method to calculate $c_{0}$.

The matching condition is that the near-field behavior of the outer expansion as $(r,\theta_\f^{i} ,\varphi_\f^{i} ,\theta_\m^{j} ,\varphi_\m^{j} )\to(R,\theta_{0},\varphi_{0},\theta_{0},\varphi_{0})$ must agree with the far-field behavior of the inner expansion as $\rho\to\infty$. That is,
\begin{align*}
\eps^{-3}v_{0}+v_{1}+\cdots
&\sim\eps^{-3}w_{0}+\cdots,\;
\text{as }(r,\theta_\f^{i} ,\varphi_\f^{i} ,\theta_\m^{j} ,\varphi_\m^{j} )\to(R,\theta_{0},\varphi_{0},\theta_{0},\varphi_{0}),\;\rho\to\infty.
\end{align*}
Using \eqref{innersol0} and \eqref{w0far}, it follows that
\begin{align*}
\alpha = v_0,
\end{align*}
and that $v_{1}$ has the singular behavior as $(r,\theta_\f^{i} ,\varphi_\f^{i} ,\theta_\m^{j} ,\varphi_\m^{j} )\to(R,\theta_{0},\varphi_{0},\theta_{0},\varphi_{0})$,
\begin{align}
\begin{split}\label{sing}
v_{1}
\sim
-v_{0}c_{0}\Big[(\tfrac{r}{R}-1)^{2}
&+c_{11}^2(\theta_\f^{i} -\theta_{0})^{2}+c_{11}^2\sin^{2}(\theta_{0})(\varphi_\f^{i} -\varphi_{0})^{2}\\
&+c_{22}^2 \big( (\theta_\m^{j} -\theta_{0}) + \frac{c_{21}}{c_{22}}(\theta_\f^{i} - \theta_0) \big)^{2}\\
& +c_{22}^2 \sin^2(\theta_0) \big( (\varphi_\m^{j} -\varphi_{0}) + \frac{c_{21}}{c_{22}}(\varphi_\f^{i} - \varphi_0) \big)^{2} \Big]^{-3/2}.
\end{split}
\end{align}


\subsection{Distributional form of the singularity}

Writing the singular behavior \eqref{sing} in distributional form for each $i\in\{1,\dots,N_{A}\}$ and $j\in\{1,\dots,N_{B}\}$, the problem in \eqref{v1first} becomes
\begin{align}
0
&=(\Dtr\Delta_{0}+\Deff_\f \L_\f+\Deff_\m \L_\m)v_{1},\quad r>R,\label{v1pde}\\
\partial_{r}v_{1}
&=\frac{v_{0}K_{0}}{R}
\sum_{i=1}^{N_{A}}\sum_{j=1}^{N_{B}}
\frac{\delta(\theta_{\f}^{i}-\theta_{0})}{\sin(\theta_{0})}\delta(\varphi_{\f}^{i}-\varphi_{0})
\frac{\delta(\theta_{\m}^{j}-\theta_{0})}{\sin(\theta_{0})}\delta(\varphi_{\m}^{j}-\varphi_{0}),
\quad r=R,\label{v1bc}
\end{align}
where
\begin{align*}
K_{0}
:=\frac{4\pi^{2}c_{0}}{ c_{11}^2 c_{22}^2}.
\end{align*}

To derive the distributional form \eqref{v1bc} of the singular behavior \eqref{sing}, assume that a function $f$ satisfies \eqref{v1pde}-\eqref{v1bc}. To derive the singular behavior of $f$ as $(r,\theta_\f^{i} ,\varphi_\f^{i} ,\theta_\m^{j} ,\varphi_\m^{j} )\to(R,\theta_{0},\varphi_{0},\theta_{0},\varphi_{0})$, define the inner solution $g$ analogously to \eqref{innersol0} and introduce the inner expansion
\begin{align}\label{gexp}
g\sim\eps^{-3}g_{0}+\cdots.
\end{align}
By the same argument that led to \eqref{harmonic}, we have that $g_{0}$ is harmonic in the five variables $({{z}},{{x_{\f}}},{{y_{\f}}},{{x_{\m}}},{{y_{\m}}})$ for $z>0$,
\begin{align}\label{gharmonic}
\big(\partial_{{{z}}{{z}}}+
\partial_{{{x_{\f}}}{{x_{\f}}}}+
\partial_{{{y_{\f}}}{{y_{\f}}}}+
\partial_{{{x_{\m}}}{{x_{\m}}}}+
\partial_{{{y_{\m}}}{{y_{\m}}}}\big)g_{0}=0,\quad {{z}}>0,\,({{x_{\f}}},{{y_{\f}}},{{x_{\m}}},{{y_{\m}}})\in\R^{4}.
\end{align}
Furthermore, $g_{0}$ satisfies the following boundary condition at ${{z}}=0$,
\begin{align}\label{gbc}
\partial_{z}g_{0}=v_{0}K_{0}c_{11}^2 c_{22} ^2 \delta(x_{A})\delta(y_{A})\delta(x_{B})\delta(y_{B}).
\end{align}
To derive \eqref{gbc}, first recall that $f$ satisfies the boundary condition in \eqref{v1bc}, then use that $g$ is defined analogously to $w$ in \eqref{innersol0}, and finally use the expansion in \eqref{gexp} to obtain
\begin{align*}
\partial_{z} g_0
&= \frac{v_0 K_0 \eps^{4}}{\sin^2(\theta_0)} \delta \Big( \eps \frac{x_{A}}{c_{11}} \Big) \delta \Big( \eps \frac{y_{A}}{c_{11} 
\sin(\theta_0)} \Big) \delta \Big( \eps \frac{ (c_{11} x_{B} - c_{21} x_{A})}{c_{11} c_{22} } \Big) \delta \Big( \eps \frac{ (c_{11} y_{B} - c_{21} y_{A})}{c_{11} c_{22} \sin(\theta_0)} \Big).
\end{align*}
Using the identity $\delta(\alpha x)=\delta(x)/|\alpha|$ and simplifying then yields \eqref{gbc}.

The solution to \eqref{gharmonic}-\eqref{gbc} is 
\begin{align}\label{g0soln}
g_{0}=\frac{-1}{4\pi^{2}}\left(v_{0}K_{0} c_{11}^2 c_{22}^2 \right)(z^{2}+x_{A}^{2}+y_{A}^{2}+x_{B}^{2}+y_{B}^{2})^{-3/2}.
\end{align}
Matching the far-field behavior of $\eps^{-3}g_{0}$ with the near-field behavior of $f$ shows that $f$ indeed has the singular behavior in \eqref{sing}. To derive \eqref{g0soln}, note that the 5D Laplacian Green's function is
\begin{align*}
G(\mathbf{a},\mathbf{b})
&=G\big((a_{1},a_{2},a_{3},a_{4},a_{5}),
(b_{1},b_{2},b_{3},b_{4},b_{5})\big)\\
&:=\frac{-1}{8\pi^{2}}\Big(\sum_{i=1}^{5}(a_{i}-b_{i})^{2}\Big)^{-3/2}
=\frac{-1}{8\pi^{2}\|\mathbf{a}-\mathbf{b}\|^{3}},
\end{align*}
and satisfies $\Delta G = \delta(\mathbf{a}-\mathbf{b})$. To derive this, we integrate over a 5D sphere centered at $\mathbf{b}$ (denoted by $B(\mathbf{b})$), to obtain
\begin{align*}
\int_{B(\mathbf{b})} \Delta G\,\dd \mathbf{a}
&= \int_{\partial B(\mathbf{b})} \partial_{\mathbf{n}} G\,\dd S_{\mathbf{a}}
=\frac{3}{8\pi^{2}} \int_{\partial B(\mathbf{b})} \frac{1}{\|\mathbf{a}-\mathbf{b}\|^{4}}\,\dd S_{\mathbf{a}}
=1,
\end{align*}
where $\dd S_{\mathbf{a}}$ denotes the surface element and we have used that the surface area of a 5D unit sphere is $\frac{8}{3}\pi^{2}$. The solution \eqref{g0soln} follows.

\subsection{Finding the capacitance $C$ and bimolecular reaction rate $k_{0}$}

Integrating the PDE \eqref{v1pde} over the region
\begin{align*}
(r,\theta_\f ,\varphi_\f ,\theta_\m ,\varphi_\m ,\theta_{0},\varphi_{0})
\in (R,R')\times \Sigma,\quad\text{where }\Sigma:=\big([0,\pi]\times[0,2\pi)\big)^{1+N_{A}+N_{B}},
\end{align*}
for $R'>R$ and using the divergence theorem and the boundary condition \eqref{v1bc}, we obtain
\begin{align}
\begin{split}\label{big0}
0
&=\int_{\Sigma}\int_{R}^{R'}(\Dtr\Delta_{0}+\Deff_\f \LL_\f+\Deff_\m \LL_\m)v_{1}\,r^{2}\dd r\,\dd\Sigma\\
&={\Dtr}\Big[\int_{\Sigma}\partial_{r}v_{1}|_{r=R'}\,R'^{2}\,\dd\Sigma-\int_{\Sigma}\partial_{r}v_{1}|_{r=R}\,R^{2}\,\dd\Sigma\Big]\\
&={\Dtr}\Big[\int_{\Sigma}\partial_{r}v_{1}|_{r=R'}\,R'^{2}\,\dd\Sigma-(4\pi)^{N_{A}+N_{B}-1} v_{0}K_{0}RN_{A}N_{B} \Big],
\end{split}
\end{align}
where
\begin{align*}
\dd\Sigma
:=\Big(\sin\theta_{0}\,\dd\theta_{0}\,\dd\varphi_{0}\Big)\Big(\prod_{i=1}^{N_{A}}\sin\theta_{A}^{i}\,\dd\theta_{A}^{i}\,\dd\varphi_{A}^{i}\Big)\Big(\prod_{j=1}^{N_{B}}\sin\theta_{B}^{j}\,\dd\theta_{B}^{j}\,\dd\varphi_{B}^{j}\Big).
\end{align*}

Now, by the far-field behavior of $p$ in \eqref{pffc}, the definition of $v$ in \eqref{vdef}, and the expansion in \eqref{outer}, it follows that $v_{1}\sim -1/r$ as $r\to\infty$. Hence, 
\begin{align*}
\int_{\varphi_{0},\varphi_\f ,\varphi_\m =0}^{2\pi}\int_{\theta_{0},\theta_\f ,\theta_\m =0}^{\pi}\partial_{r}v_{1}|_{r=R'}\,R'^{2}\,\dd\Sigma
\to (4\pi)^{N_{A}+N_{B}+1},\quad\text{as }R'\to\infty.
\end{align*}
Therefore, taking $R\to\infty$ in \eqref{big0} and solving for $v_{0}$ yields
\begin{align*}
v_{0}
=\frac{-(4\pi)^{2}}{N_\f N_\m  RK_{0}}
=\frac{-(4\pi)^{2} c_{11}^2 c_{22}^2 }{4 N_\f N_\m  R \pi^{2}c_{0}}
=\frac{-4c_{11}^2 c_{22}^2}{N_\f N_\m  Rc_{0}}.
\end{align*}
Therefore, \eqref{pffc}-\eqref{outer} yields the leading order behavior of $C$ as $\eps\to0$,
\begin{align*}
C
\sim\eps^{3}\frac{N_{A}N_{B}Rc_{0}}{4 c_{11}^2 c_{22}^2 }
=\eps^{3}N_{A}N_{B}\chi R,\quad\text{as }\eps\to0,
\end{align*}
where we have defined
\begin{align}\label{chi}
\chi
:=\frac{c_{0}}{4c_{11}^{2}c_{22}^{2}}
=(D_{A}D_{B}+D_{A}+D_{B})\frac{c_{0}}{4}.
\end{align}
Finally, upon using the relation in \eqref{k0C}, we obtain the asymptotic behavior of the bimolecular reaction rate constant,
\begin{align}\label{k0calc}
k_{0}
\sim\eps^{3}N_{A}N_{B}\chi\ks
,\quad\text{as }\eps\to0.
\end{align}
{\black Note that $\chi>0$ is a dimensionless constant which measures the how the ratios of the diffusivities $D_{A}$ and $D_{B}$ (see \eqref{Ddef}) and the relative binding site sizes $a_{A}$ and $a_{B}$ (see \eqref{radii}) affect the bimolecular binding rate in the limit of small binding sites.}

The calculation above was for the special case $\Drot_{A}=\Drot_{B}=0$ with
\begin{align}\label{Deffg0}
\Deff_{A}
&:=\Drot_{A}+R^{-2}\Dsurf_{A}>0,\quad
\Deff_{B}
:=\Drot_{B}+R^{-2}\Dsurf_{B}>0.
\end{align}
{\black However}, the final result \eqref{k0calc} still holds in the general case that
\begin{align}\label{Dg}
\Drot_{A}\ge0, \,\Drot_{B}\ge0,\, \Dsurf_{A}\ge0,\,\Dsurf_{B}\ge0,
\end{align}
as long as \eqref{Deffg0} holds. 

To see why this is the case, note that including rotational diffusion merely introduces correlations in the SDEs in \eqref{ASDE} and \eqref{BSDE}. That is, if $\Drot_{A}>0$, then the position of the $i$th $A$ binding site, $(\Theta_{A}^{i}(t),\Phi_{A}^{i}(t))$, and the position of the $j$th $A$ binding site, $(\Theta_{A}^{j}(t),\Phi_{A}^{j}(t))$, are no longer independent (since their positions depend on the rotational path of the $A$ molecule, which is common to both binding sites). These correlations in binding site positions would change the PDE satisfied by $p$ in \eqref{ppde}. However, our analysis above shows that the leading order result in \eqref{k0calc} is independent of the arrangement of binding sites (as long as they are well-seperated), and therefore \eqref{k0calc} must still hold in the general case of \eqref{Deffg0}-\eqref{Dg}.

\section{A kinetic Monte Carlo method for calculating $\chi$} \label{calculatingc0}

The asymptotic behavior of the bimolecular binding rate constant $k_{0}$ given in \eqref{k0calc} depends on $\chi$, which depends on the constant $c_0$ in the far-field behavior \eqref{w0far} of the leading order inner solution $w_{0}$ satisfying \eqref{harmonic}-\eqref{w0bc}. In particular, $c_{0}$ is the electrostatic capacitance of the 4D region $\mathcal{R}$ in \eqref{region} embedded in $\R^{5}$. Notice that $c_{0}=c_{0}(D_{A},D_{B},a_{B})$ is a function of the following three dimensionless parameters
\begin{align*}
D_{A}:=R^{2}\Deff_\f/\Dtr>0,
\quad
D_{B}:=R^{2}\Deff_\m/\Dtr>0,
\quad a_{B}\in(0,1],
\end{align*}
since we can without loss of generality take $a_{B}\le a_{A}=1$.

In this section, we develop a kinetic Monte Carlo method for rapid numerical calculation of $c_{0}$. Our approach uses a recent algorithm that was devised by Bernoff, Lindsay, and Schmidt to calculate the capacitance of 2D regions embedded in $\R^{3}$ \cite{bernoff2018b}.

\subsection{Probabilistic interpretation}\label{probint}

The method relies on a probabilistic interpretation of the PDE boundary value problem \eqref{harmonic}-\eqref{w0bc} satisfied by the leading order inner solution $w_{0}$. Let
\begin{align}\label{ZZ}
\Z(t)
=\big(Z(t),X_{\f}(t),{{Y_{\f}}}(t),{{X_{\m}}}(t),{{Y_{\m}}}(t)\big)\in\R^{5}
\end{align}
be a standard 5D Brownian motion. Define the first time that this process reaches the region $\mathcal{R}$ in \eqref{region},
\begin{align}\label{tau0}
\tau_{0}
:=\inf\{t>0:\Z(t)\in\RR\}.
\end{align}
It is straightforward to show that the leading order inner solution $w_{0}$ satisfying \eqref{harmonic}-\eqref{w0bc} can be written as
\begin{align}\label{qw}
w_{0}(z,x_{\f},y_{\f},x_{\m},{{y_{\m}}})
=v_{0}(1-q(z,x_{\f},y_{\f},x_{\m},{{y_{\m}}})),
\end{align}
where $q$ is the probability that $\Z$ eventually reaches $\RR$, conditioned on the initial position of $\Z$,
\begin{align*}
q(z,x_{\f},y_{\f},x_{\m},{{y_{\m}}})
:=\P\big(\tau_{0}<\infty\,|\,\Z(0)=(z,x_{\f},y_{\f},x_{\m},{{y_{\m}}})\big).
\end{align*}
Furthermore, the function $q$ must be harmonic for $z\neq0$, which in 5D spherical coordinates is
\begin{align}\label{qharmonic}
\Big(\frac{4}{\rho}\partial_{\rho}+\partial_{\rho\rho}+\frac{1}{\rho^{2}}\L^{(4)}\Big)q=0,\quad z\neq0,
\end{align}
where $\rho:=\sqrt{z^{2}+x_{\f}^{2}+y_{\f}^{2}+x_{\m}^{2}+{{y_{\m}}}^{2}}$ is the radius and $\L^{(4)}$ denotes the Laplace-Beltrami operator on the 4-sphere. Further, $q$ satisfies the boundary conditions at $z=0$, 
\begin{align}\label{qbc}
\begin{split}
q
&=1,\quad {{z}}=0,\; {{x_{\f}}}^{2}+{{y_{\f}}}^{2}<c_{11}^{2} a_\f^2,\; \Big({{x_{\m}}}-\frac{c_{21}}{c_{11}}{{x_{\f}}}\Big)^{2}+\Big({{y_{\m}}}-\frac{c_{21}}{c_{11}}{{y_{\f}}}\Big)^{2}<c_{22}^{2} a_\m ^2\\
\partial_{{{z}}}q
&=0,\quad {{z}}=0,\; \text{otherwise}.
\end{split}
\end{align}

Let $\overline{q}(\rho)$ denote the average of $q$ over the surface of the 5D ball of radius $\rho>0$ centered at the origin. Now, notice that if $z=0$ and $\rho>0$ is such that
\begin{align*}
\rho
>\rho_{0}:=c_{11}^{2} a_\f^2+\Big(\Big(\frac{c_{21}}{c_{11}}\Big)c_{11}a_{\f}+c_{22}a_{\m}\Big)^{2},
\end{align*}
then \eqref{qbc} ensures $\partial_{z}q=0$. In particular, $\rho_{0}$ is the smallest radius which guarantees the reflecting boundary condition in \eqref{qbc} is satisfied. Therefore, integrating \eqref{qharmonic} over the surface of the 5D ball of radius $\rho>\rho_{0}$ centered at the origin, using the divergence theorem, and interchanging integration with differentiation yields the following ODE for $\overline{q}(\rho)$,
\begin{align}\label{qode}
\Big(\frac{4}{\rho}\partial_{\rho}+\partial_{\rho\rho}\Big)\overline{q}=0,\quad \rho>\rho_{0}.
\end{align}
The general solution to \eqref{qode} is $\overline{q}(\rho)=K_{1}\rho^{-3}+K_{2}$ for constants $K_{1},K_{2}\in\R$. The relation \eqref{qw} and the far-field behavior of $w_{0}$ in \eqref{w0far} implies $K_{2}=0$ and $K_{1}=c_{0}$.

\subsection{Kinetic Monte Carlo algorithm}

We have shown in the previous subsection that we can find $c_{0}$ by calculating the probability, $\overline{q}(\rho)$, that the 5D Brownian motion $\Z$ in \eqref{ZZ} eventually reaches the region $\mathcal{R}$ defined by \eqref{region}, conditioned that $\Z$ is initially uniformly distributed on a ball of radius $\rho>\rho_{0}$. {\black Roughly speaking, we therefore} approximate $\overline{q}(\rho)$ by simulating $M\gg1$ diffusive paths of $\Z$ and calculating the proportion of these $M$ paths which reach $\mathcal{R}$ before some large outer radius $\rho_{\infty}\gg\rho_{0}$.

However, simulating these diffusive paths with a standard time discretization scheme would be incredibly computationally expensive. Indeed, the Brownian motion would have to take many steps to reach the outer radius $\rho_{\infty}$ unless the discrete time step $\Delta t>0$ is taken very large. On the other hand, the time step $\Delta t>0$ would need to be taken very small in order to accurately resolve the dynamics of $\Z$ near $\mathcal{R}$.

We therefore develop a kinetic Monte Carlo method which avoids these issues \cite{bernoff2018b}. This kinetic Monte Carlo method breaks the simulation process into two steps, where each step corresponds to a simpler diffusion problem that can be exactly and efficiently simulated. The method then alternates between these steps until the simulation reaches a break point. The method takes very large time steps and generates statistically exact paths of $\Z$. Indeed, in calculating $c_{0}$ from this method, the only error stems from the finite outer radius $\rho_{\infty}<\infty$ and the finite number of diffusive paths $M<\infty$ (as opposed to error stemming from a nonzero time step). Furthermore, the computational efficiency of the method allows us to mitigate these two sources of error by taking $\rho_{\infty}$ and $M$ very large. For example, simulating $M=10^{6}$ paths with $\rho_{\infty}=10^{5}$ takes roughly 10 seconds on a standard personal laptop computer.

To describe the method, notice that the 5D Brownian motion $\Z$ in \eqref{ZZ} can be visualized as a pair of 3D Brownian motions,
\begin{align*}
\X_{\f}(t)
&:=(X_{\f}(t),{{Y_{\f}}}(t),Z(t))\in\R^{3},\\
{{\X_{\m}}}(t)
&:=({{X_{\m}}}(t),{{Y_{\m}}}(t),Z(t))\in\R^{3},
\end{align*}
with independent $x$ and $y$ coordinates and identical $z$ coordinates. Therefore, the 5D Brownian motion in $\Z$ reaches the region $\mathcal{R}$ in \eqref{region} if and only if $\X_{\f}$ and ${{\X_{\m}}}$ reach the $z=0$ plane in $\R^{3}$ while (i) $\X_{\f}$ is in a disk of radius $c_{11}a_{\f}$ centered at the origin and (ii) ${{\X_{\m}}}$ is in a disk of radius $c_{22}a_{\m}$ centered at $\frac{c_{21}}{c_{11}}\X_{\f}$.

After initially placing the ``particle'' $\Z$ on the 5D sphere of radius $\rho$ centered at the origin according to a uniform distribution, the method employs the following two stages developed by Bernoff, Lindsay, and Schmidt \cite{bernoff2018b} (originally developed for 3D diffusion). {\black We note that Stage II is the classical ``walk-on-spheres'' method due to Muller in 1956 \cite{muller1956}.}
\begin{itemize}
\item
\emph{Stage \textup{I:} Projection from bulk to plane.} The particle is projected to the $z=0$ plane following the exact distribution given below. If the particle lands in $\RR$, then this event is recorded and the trial ends. If not, the algorithm proceeds to Stage II. 

\item
\emph{Stage \textup{II:} Projection from plane to the bulk.} A distance $\nu>0$ is calculated which is less than or equal to the distance from the current particle location to $\mathcal{R}$. The particle is then projected to a uniformly chosen random point on the 5D sphere with radius $\nu$, centered at the current particle location. If the particle reaches a distance that is larger than $\rho_{\infty}$ from the origin, then this event is recorded and the trial ends. If not, the algorithm returns to Stage I.
\end{itemize}

{\black 
We now describe the basic idea behind the method. The method aims to simulate whether a diffusing particle eventually reaches the region $\mathcal{R}$. In order to reduce computational cost, the method skips
directly simulating intermediate steps of the diffusing particle before the particle could possibly reach $\mathcal{R}$. Since $\mathcal{R}$ is a subset of the $z=0$ plane, Stage I skips simulating all the steps until the particle reaches $z=0$. Then, if the particle is in $\mathcal{R}$, the simulation ends. If the particle is not in $\mathcal{R}$, then it is some positive distance $\nu_{0}>0$ away from $\mathcal{R}$. So, Stage II moves the particle a distance $\nu$ with $\nu\in(0,\nu_{0}]$ (calculating the exact distance $\nu_{0}$ is difficult). With probability one, this new position is not on the $z=0$ plane and is thus in the ``bulk.'' At this point, if the particle is a large distance ($\rho_{\infty}$) from the origin, then we assume that the particle will never reach $\mathcal{R}$, and so the simulation ends. Otherwise, the simulation returns to Stage I.
}

To calculate the distribution in Stage I, we first sample the random time it takes $\Z$ to reach $z=0$, which is \cite{bernoff2018b}
\begin{align*}
t^* = \frac{1}{4} \Big( \frac{z}{\text{erfc}^{-1}(U) } \Big)^2,
\end{align*}
where $U\sim\text{uniform}(0,1)$ is uniformly distributed on $[0,1]$. Then, if $\Z$ is at position $(z,x_{\f},y_{\f},x_{\m},{{y_{\m}}})\in\R^{5}$ at the start of Stage I, the position at the end of Stage I is
\begin{align*}
(0,x_{\f},y_{\f},x_{\m},{{y_{\m}}})+\sqrt{2t^{*}}(0,\xi_{1},\xi_{2},\xi_{3},\xi_{4})\in\R^{5}
\end{align*}
where $\xi_{1},\xi_{2},\xi_{3},\xi_{4}$ are independent standard normal random variables.

In Stage II, we want to propagate the particle as far as possible, while ensuring that the particle cannot reach $\mathcal{R}$ during this propagation \cite{bernoff2018b, muller1956}. Let $t_{0}>0$ be the time at the start of Stage II. If the algorithm is in Stage II, then it must be the case that $\Z(t_{0})\notin\RR$, and thus
\begin{align*}
d_{1}:=\|\X_{A}(t_{0})\|-r_{1}>0,\quad\text{and/or}\quad
d_{2}:=\|\X_{B}(t_{0})-s\X_{A}(t_{0})\|-r_{2}>0,
\end{align*}
where
\begin{align*}
r_{1}:=c_{11}a_{\f},\quad
r_{2}:=c_{22}a_{\m},\quad
s:=c_{21}/c_{11},
\end{align*}
and $\|\cdot\|$ denotes the standard Euclidean norm. Now, if $\Z(t_{1})\in\RR$ for $t_{1}>t_{0}$ and $d_{1}>0$, then it must be the case that
\begin{align*}
\|\X_{A}(t_{1})-\X_{A}(t_{0})\|\ge d_{1}.
\end{align*}
Similarly, if $\Z(t_{1})\in\RR$ for $t_{1}>t_{0}$ and $d_{2}>0$, then it must be the case that
\begin{align}\label{constraint}
s\|\X_{A}(t_{1})-\X_{A}(t_{0})\|+\|\X_{B}(t_{1})-\X_{B}(t_{0})\|\ge d_{2}.
\end{align}
A straightforward calculus exercise shows that the minimum distance, $\|\Z(t_{1})-\Z(t_{0})\|$, subject to the constraint \eqref{constraint} is
\begin{align*}
\frac{d_{2}}{\sqrt{1+s^{2}}}>0.
\end{align*}
Therefore, if we define the distance
\begin{align*}
\nu
:=\max\Big\{d_{1},\frac{d_{2}}{\sqrt{1+s^{2}}}\Big\},
\end{align*}
then {\black it} follows that the 5D sphere of radius $\nu$ centered at $\Z(t_{0})\notin\RR$ cannot intersect $\RR$. Hence, Stage II places the particle uniformly on the boundary of this 5D sphere (the uniform distribution follows from symmetry of Brownian motion).

\begin{figure}
  \centering
    \includegraphics[width=	1\textwidth]{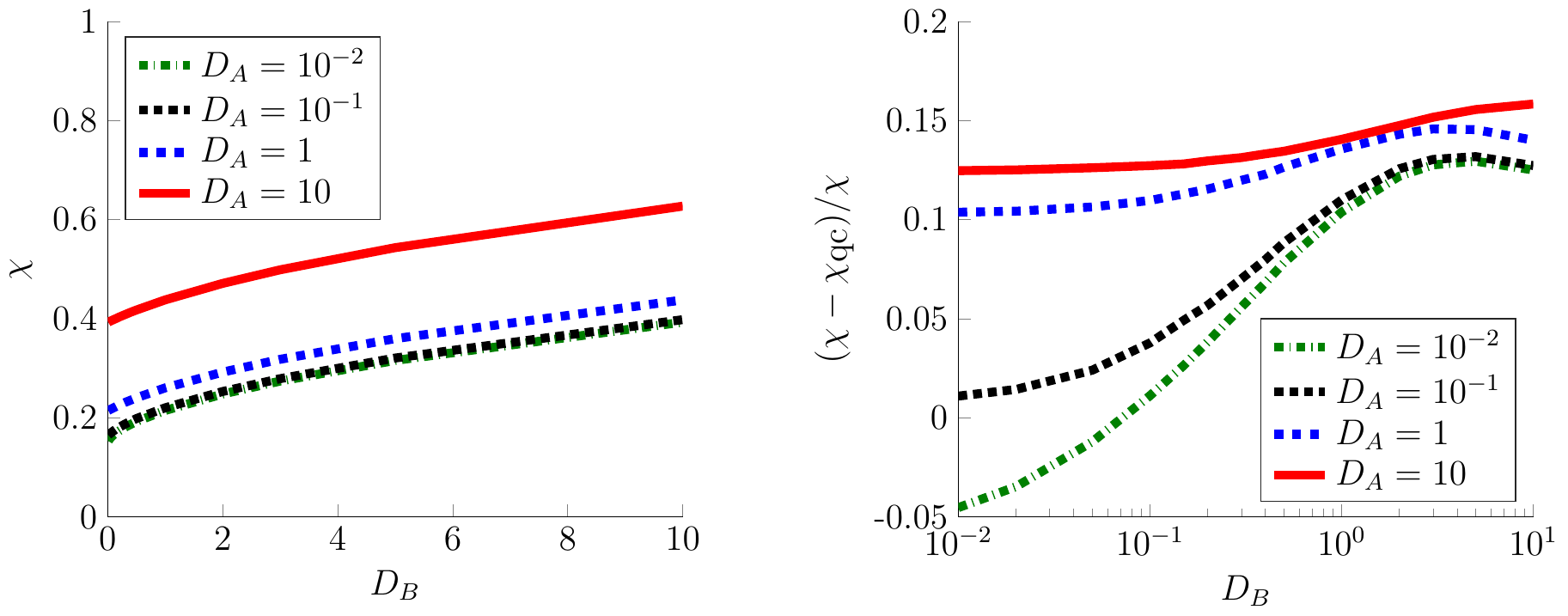}
 \caption{\textbf{Left:} The factor $\chi>0$ as a function of $D_{B}:=R^{2}\Deff_{B}/\Dtr$ for different values of $D_{A}:=R^{2}\Deff_{A}/\Dtr$. The curves for $D_{A}=10^{-2}$ and $D_{A}=10^{-1}$ are almost indistinguishable. \textbf{Right:} The relative error between $\chi$ and the approximation $\chiqc$ in \eqref{chia2}. In both panels, the values of $\chi$ are computed from $M=10^{8}$ trials of the kinetic Monte Carlo method of section~\ref{calculatingc0} with outer radius $\rho_{\infty}=10^{5}$ and $a_{A}=a_{B}=1$.}
 \label{figchi}
\end{figure}

In Figure~\ref{figchi}, we plot $\chi=(D_{A}D_{B}+D_{A}+D_{B})c_{0}/4$ (see \eqref{chi}) as a function of $D_{A}$ and $D_{B}$ using the above kinetic Monte Carlo method. For each pair of $D_{A},D_{B}$, the value of $c_{0}$ used in $\chi$ is computed from $M=10^{8}$ trials with outer radius $\rho_{\infty}=10^{5}$. This figure shows that $\chi$ is an increasing function of $D_{A},D_{B}$ (as expected) and that $\chi$ varies between roughly $\chi\approx0.17$ and $\chi\approx0.63$ for $(D_{A},D_{B})\in[10^{-2},10]^{2}$. In particular, $\chi$ varies by less than a factor of 4 as $D_{A}$ and $D_{B}$ each vary 3 orders of magnitude.

Notice that the symmetry in the full binding model of section~\ref{math} implies that $\chi$ must be symmetric in $D_{A},D_{B}$ (that is, $\chi(D_{A},D_{B})=\chi(D_{B},D_{A})$). To test this, we computed $\chi$ for
\begin{align}\label{pairs}
(D_{A},D_{B})\in\{ 0.01, 0.02, 0.05, 0.1, 0.15,0.2, 0.3, 0.4, 0.5, 1,2,3,5,10 \}^{2}\subset\R^{2}
\end{align}
for a total of $14^{2}=196$ distinct pairs of $D_{A}$ and $D_{B}$ values (we take $a_{A}=a_{B}=1$). For each pair, $c_{0}$ is computed from $M=10^{8}$ trials with outer radius $\rho_{\infty}=10^{5}$. Using this data, the maximum relative difference, $|\chi(D_{A},D_{B})-\chi(D_{B},D_{A})|/\chi(D_{A},D_{B})$, for the 196 pairs in \eqref{pairs} is 0.0013, which is well within the expected error due to $M=10^{8}<\infty$ trials (see section~\ref{accuracy} below). This symmetry is a necessary self-consistency check, as it is not \emph{a priori} clear from the PDE boundary value problem in \eqref{harmonic}-\eqref{w0bc} that $c_{0}$ is symmetric in $D_{A},D_{B}$ (though these simulations indicate that it is). 

\subsection{Accuracy}\label{accuracy}

In calculating $c_{0}$ from the method described above, the only error stems from the finite outer radius $\rho_{\infty}<\infty$ and the finite number of diffusive paths $M<\infty$. In this subsection, we estimate the error as a function of $\rho_{\infty}$ and $M$.

If $\P_{\rho}$ denotes the probability measure conditioned that $\Z$ starts uniformly on the 5D sphere of radius $\rho>\rho_{0}$ centered at the origin, then it follows from the analysis in section~\ref{probint} that
\begin{align}
\begin{split}\label{infeq}
\overline{q}(\rho)
=\frac{c_{0}}{\rho^{3}}
=\P_{\rho}(\tau_{0}<\infty)
&=\P_{\rho}(\tau_{0}<\tau_{\rho_{\infty}}<\infty)+\P_{\rho}(\tau_{\rho_{\infty}}<\tau_{0}<\infty)\\
&=\P_{\rho}(\tau_{0}<\tau_{\rho_{\infty}})+\P_{\rho}(\tau_{\rho_{\infty}}<\tau_{0}<\infty),
\end{split}
\end{align}
where {\black $\tau_{0}$ is the first time $\Z$ reaches $\mathcal{R}$ (see \eqref{tau0}) and} $\tau_{\rho_{\infty}}$ is the first time $\Z$ reaches distance $\rho_{\infty}$ from the origin (the final equality in \eqref{infeq} follows from the fact that $\tau_{\rho_{\infty}}<\infty$ with probability one).

Notice that the numerical algorithm actually approximates {\black a probability that is bounded above by $\P_{\rho}(\tau_{0}<\infty)$ and below by $\P_{\rho}(\tau_{0}<\tau_{\rho_{\infty}})$. In particular, it is bounded above by $\P_{\rho}(\tau_{0}<\infty)$ since the algorithm neglects some paths which first reach distance $\rho_{\infty}$ and then reach $\mathcal{R}$ (if the algorithm terminates in Stage II). Further, it is bounded below by $\P_{\rho}(\tau_{0}<\tau_{\rho_{\infty}})$ since the algorithm includes some paths which first reach distance $\rho_{\infty}$ and then reach $\mathcal{R}$ (since the particle may reach distance $\rho_{\infty}$ in Stage I before reaching the $z=0$ plane).} In light of \eqref{infeq}, we thus want to show that $\P_{\rho}(\tau_{\rho_{\infty}}<\tau_{0}<\infty)$ is small if $\rho_{\infty}\gg\rho$. Notice that $\P_{\rho}(\tau_{\rho_{\infty}}<\tau_{0}<\infty)$ is the probability of paths that first reach distance $\rho_{\infty}$ from the origin, and then {\black reach} $\mathcal{R}$. Since particles that start at distance $\rho_{\infty}$ from the origin have probability roughly $c_{0}\rho_{\infty}^{-3}$ of reaching $\mathcal{R}$, it follows that $\P_{\rho}(\tau_{\rho_{\infty}}<\tau_{0}<\infty)$ decays like $\rho_{\infty}^{-3}$ as $\rho_{\infty}$ grows.

More precisely, subtracting $\P_{\rho}(\tau_{0}<\tau_{\rho_{\infty}})$ from \eqref{infeq}, multiplying by $\rho^{3}$ and dividing by $c_{0}$ yields
\begin{align}\label{bounds}
0
\le\frac{c_{0}-\rho^{3}\P_{\rho}(\tau_{0}<\tau_{\rho_{\infty}})}{c_{0}}
=\frac{\rho^{3}\P_{\rho}(\tau_{\rho_{\infty}}<\tau_{0}<\infty)}{c_{0}}.
\end{align}
Now, it follows from the strong Markov property that
\begin{align*}
\P_{\rho}(\tau_{\rho_{\infty}}<\tau_{0}<\infty)
\le\inf_{\x\in\R^{5}:\|\x\|=\rho_{\infty}}\P(\tau_{0}<\infty\,|\,\Z(0)=\x)
=\frac{c_{0}}{\rho_{\infty}^{3}}+o(\rho_{\infty}^{-3}),\,\text{as }\rho_{\infty}\to\infty.
\end{align*}
Therefore, \eqref{bounds} implies that the relative error between $c_{0}$ {\black and} $\rho^{3}\P_{\rho}(\tau_{0}<\tau_{\infty})$ decays like $\rho_{\infty}^{-3}$,
\begin{align*}
0
\le\frac{c_{0}-\rho^{3}\P_{\rho}(\tau_{0}<\tau_{\rho_{\infty}})}{c_{0}}
=\Big(\frac{\rho}{\rho_{\infty}}\Big)^{3}+o(\rho_{\infty}^{-3}),\quad\text{as }\rho_{\infty}\to\infty.
\end{align*}
In our simulations, we take $\rho_{\infty}=10^{5}$, and $\rho$ to be order one, which means the relative error in approximating $c_{0}$ that stems from $\rho_{\infty}<\infty$ is on the order of $10^{-15}$. As an aside, we note that we could take smaller values of $\rho_{\infty}$ and still obtain very accurate results, but the computational cost of the algorithm depends very weakly on $\rho_{\infty}$. 

Estimating the error in the approximation that stems from a finite number of trials, $M<\infty$, is a basic problem in statistical inference. Each kinetic Monte Carlo trial is an independent realization of a Bernoulli random variable with parameter $p_{0}\textcolor{black}{\in(\P_{\rho}(\tau_{0}<\tau_{\rho_{\infty}}),\P_{\rho}(\tau_{0}<\infty))\subset}(0,1)$, and we are estimating $p_{0}$ by the fraction $p_{\text{kmc}}\in[0,1]$ of trials {\black which terminate in Stage I}. Given an estimate $p_{\text{kmc}}$ formed from $M$ trials, the $100(1-\alpha)\%$ confidence interval for $p_{0}$ can be estimated by \cite{agresti1998}
\begin{align*}
p_{\pm}
:=\left(p_{\text{kmc}}
+\frac{z_{\alpha/2}^{2}}{2M}
\pm z_{\alpha/2}\sqrt{\frac{p_{\text{kmc}}(1-p_{\text{kmc}})+z_{\alpha/2}^{2}/(4M)}{M}}\right)\left(1+\frac{z_{\alpha/2}^{2}}{M}\right)^{-1},
\end{align*}
where $z_{c}$ denotes the $1-c$ quantile of the standard normal distribution. That is, $p_{0}\in[p_{-},p_{+}]$ with approximate probability $1-\alpha$. Applying this statistical test to our simulations which use $M=10^{8}$ trials, we find that for each of the 196 choices of parameters $(D_{A},D_{B})$ in \eqref{pairs}, our estimate of $c_{0}$ has a relative error of less than $0.002$ with probability $1-\alpha=0.95$.

\subsection{Stokes-Einstein relation}

In this subsection, we briefly discuss how $\chi$ varies as a function of the relative sizes of the $A$ and $B$ molecules if we assume that the Stokes-Einstein relation holds and that there is no surface diffusion $\Dsurf_{A}=\Dsurf_{B}=0$. In particular, the Stokes-Einstein relation implies
\begin{align*}
\Dtr
=\Dtr_{A}+\Dtr_{B}
=\frac{k_{\text{B}}T}{6\pi\eta }\Big(\frac{1}{R_{A}}+\frac{1}{R_{B}}\Big),
\quad
\Drot_{A}
=\frac{k_{\text{B}}T}{8\pi\eta (R_{A})^{3}},
\quad
\Drot_{B}
=\frac{k_{\text{B}}T}{8\pi\eta (R_{B})^{3}},
\end{align*}
where $k_{\text{B}}$ is Boltzmann's constant, $T$ is temperature, and $\eta$ is the viscosity of the medium. Therefore, recalling that $R:=R_{A}+R_{B}$, we obtain that $D_{A}$ and $D_{B}$ are merely geometric factors \cite{Berg1985, lawley2019bp},
\begin{align*}
D_{A}
=\frac{R^{2}\Drot_{A}}{\Dtr}=\frac{3}{4}{\xi}(1+{\xi}),\quad
D_{B}
=\frac{R^{2}\Drot_{B}}{\Dtr}=\frac{3}{4}{\xi}^{-1}(1+{\xi}^{-1}),
\end{align*}
where we can without loss of generality take the $B$ molecule to be smaller than the $A$ molecule,
\begin{align}\label{rho}
{\xi}:=\frac{R_{B}}{R_{A}}\le1.
\end{align}

If we further assume that the $A$ and $B$ binding sites have respective radii $\eps R_{A}$ and $\eps R_{B}$ for some $\eps\ll1$ (meaning $a_{A}=a_{B}$), then the bimolecular binding rate $k_{0}$ in \eqref{k0calc} simplifies to
\begin{align*}
k_{0}
\sim\eps^{3}N_{A}N_{B}\chi({\xi})\ks,\quad\text{as }\eps\to0,
\end{align*}
where $\chi:(0,1]\to(0,\infty)$ is a function of the single parameter ${\xi}$ in \eqref{rho}. Using the kinetic Monte Carlo method above, we find that $\chi=\chi({\xi})$ is a decreasing function of ${\xi}\in(0,1]$ and that
\begin{align}\label{chivals}
\chi(\tfrac{1}{10})\approx0.89,\quad
\chi(\tfrac{1}{4})\approx0.46,\quad
\chi(\tfrac{1}{2})\approx0.33,\quad
\chi(\tfrac{3}{4})\approx0.30,\quad
\chi(1)\approx0.29.
\end{align}
Hence, \eqref{chivals} reveals that $\chi$ varies very little as a function of ${\xi}\in(0,1]$, unless ${\xi}$ is very small.

\section{Incorporating binding site competition}\label{comp}

The asymptotic behavior of the bimolecular binding rate $k_{0}$ in \eqref{k0calc} is in the $\eps\to0$ limit. In particular, formula \eqref{k0calc} is not valid if we fix $\eps>0$ and take $N_{A}\to\infty$ and/or $N_{B}\to\infty$. Indeed, formula \eqref{k0calc} grows without bound as $N_{A}$ and/or $N_{B}$ grows, whereas $k_{0}$ must always be bounded above by $\ks$. However, it is immediately clear that if $N_{A}\to\infty$ (or $N_{B}\to\infty$), then $k_{0}$ should simply approach the binding rate for the case of one molecule completely covered by binding sites and one molecule partially covered (i.e.\ one homogeneous molecule and one heterogeneous molecule). Specifically, we expect that
\begin{align}\label{limitbp0}
k_{0}
\to k_{A}\quad\text{as }N_{B}\to\infty,
\end{align}
where $k_{A}$ is binding rate constant for a homogeneous $B$ molecule and a heterogeneous $A$ molecule. In the case of small, well-separated binding sites, this $k_{A}$ was recently shown in \cite{lawley2019bp} to be well-approximated by
\begin{align}\label{recent}
k_{A}
\approx
\overline{k_{A}}:=\frac{\lambda_{A}N_{A}a_{A}\eps}{\pi+\lambda_{A}N_{A}a_{A}\eps}\ks,\quad\text{where }\lambda_{A}:=\sqrt{1+\frac{R^{2}\Drot_{A}+\Dsurf_{A}}{\Dtr}}.
\end{align}
In fact, the limiting behavior in \eqref{recent} must also hold if $\lambda_{B}\to\infty$, since the $B$ binding sites effectively cover the $B$ molecule in this limit (see \cite{lawley2019bp,lawley2019dtmfpt,lawley2019fpk} for more on this phenomenon).

The basic reason that the asymptotic behavior in \eqref{k0calc} breaks down for fixed $\eps>0$ and sufficiently large $N_{B}$ (or sufficiently large $\lambda_{B}$, $N_{A}$, or $\lambda_{A}$) is that the binding sites begin to ``compete'' for the flux in this limit. To obtain a formula for $k_{0}$ which includes the effects of competition between binding sites, we adopt the heuristic quasi chemical formalism of {\v{S}}olc and Stockmayer's 1973 study of a single binding site model \cite{solc1973}. In addition to yielding such a formula for $k_{0}$, we find that by combining this approach with \eqref{recent}, we obtain a simple analytical approximation for $\chi$.

\subsection{Quasi chemical formalism of {\v{S}}olc and Stockmayer \cite{solc1973}}

\begin{figure}
  \centering
    \includegraphics[width=	.8\textwidth]{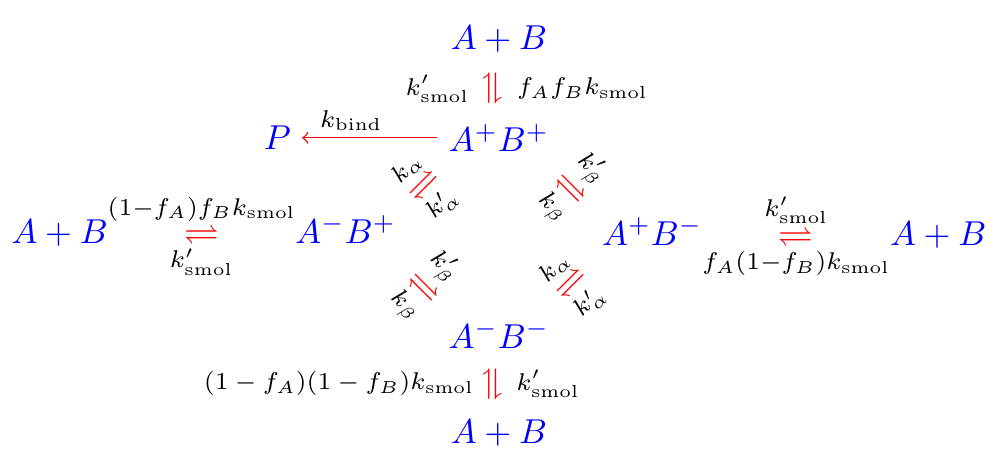}
 \caption{Chemical reaction diagram for the quasi chemical approximation.}
 \label{figstate}
\end{figure}

The quasi chemical formalism is a heuristic approximation that collapses the infinite-dimensional state space of a diffusion-based binding model into a discrete state space model with 6 states. In this discrete state model, the molecules can be \emph{far} from each other, \emph{close} to each other, or \emph{bound}, and if they are close, then we distinguish whether or not a binding site of $A$ (respectively $B$) is aligned toward $B$ (respectively $A$). This model is depicted in Figure~\ref{figstate}, where $A+B$ denotes that the molecules are far from each other, $P$ denotes that the particles have bound and formed an irreversible product, and $A^{\pm}B^{\pm}$ denotes the 4 possible close states with the $+$ superscript denoting a binding site aligned and the $-$ superscript denoting no binding site aligned. For example, $A^{+}B^{-}$ means the molecules are close and an $A$ binding site is aligned toward $B$ and no $B$ binding site is aligned toward $A$. 

The transition rates between the states are given in Figure~\ref{figstate}. Note that the rate from the far state to a close state is the Smoluchowski rate, $\ks$, multiplied by the corresponding binding site surface fractions. For example,
\begin{align*}
A+B\to A^{+}B^{-}\quad\text{at rate}\quad f_{A}(1-f_{B})\ks>0,
\end{align*}
where $f_{A}\in(0,1]$ denotes the fraction of the surface of $A$ covered in binding sites and similarly for $f_{B}\in(0,1]$. Further, $\ks'$ denotes the rate that the molecules diffuse away from each other, which is the same for the 4 states $A^{\pm}B^{\pm}$, and $\kb$ denotes the rate that the molecules bind once they are aligned ($A^{+}B^{+}\to P$). Finally, $k_{\alpha}$ (respectively $k_{\beta}$) denotes the rate that $A$ (respectively $B$) aligns toward $B$ (respectively $A$), and the reverse rate is $k_{\alpha}'=\frac{1-f_{A}}{f_{A}}k_{\alpha}$ (respectively $k_{\beta}'=\frac{1-f_{B}}{f_{B}}k_{\beta}$), which follows from microscopic reversibility \cite{solc1973}.

Writing down the system of mass action ODEs corresponding to the reaction diagram in Figure~\ref{figstate} and solving for the steady state yields the following effective bimolecular binding rate constant \cite{solc1973},
\begin{align}\label{keff0}
\keff
:=\frac{\kb [A^{+}B^{+}]}{[A][B]}
=f_{A}f_{B}\ks\big(\ks'/\kb+\Lambda_{A}\Lambda_{B}+\psi\big)^{-1},
\end{align}
where $[A][B]$ denotes the product of the steady state concentrations of $A$ and $B$ molecules which are far from each other, $[A^{+}B^{+}]$ denotes the steady state concentration $A$ and $B$ molecules which are close and aligned, and
\begin{align}\label{solc}
\begin{split}
\Lambda_{A}
&:=\frac{r_{A}+1}{r_{A}+f_{A}}f_{A},
\quad\Lambda_{B}:=\frac{r_{B}+1}{r_{B}+f_{B}}f_{B},\\
\psi
&:=\Big[\frac{1}{(1-\Lambda_{A})(1-\Lambda_{B})}
+\frac{1}{(1-\Lambda_{A})(\Lambda_{B}-f_{B})}
+\frac{1}{(1-\Lambda_{B})(\Lambda_{A}-f_{A})}\Big]^{-1},
\end{split}
\end{align}
and $r_{A}:=k_{\alpha}/\ks'$ and $r_{B}:=k_{\beta}/\ks'$. Since we assumed in previous sections that the molecules bind as soon as they are in contact, we take $\ks/\kb\to0$ in \eqref{keff0} which yields
\begin{align}\label{keff}
\keff
=\frac{f_{A}f_{B}}{\Lambda_{A}\Lambda_{B}+\psi}\ks.
\end{align}

\subsection{One homogeneous molecule and one heterogeneous molecule}

If $B$ is completely covered in binding sites (i.e.\ $f_{B}=1$), then $\Lambda_{B}=1$, $\psi=0$, and \eqref{keff} reduces to
\begin{align}\label{keffr}
\keff
=\frac{r_{A}+f_{A}}{r_{A}+1}\ks.
\end{align}
If the $N_{A}$ locally circular binding sites of radius $\eps a_{A}R_{A}>0$ are placed independently and uniformly on the surface of the $A$ molecule, then the expected surface fraction of $A$ covered in binding sites is
\begin{align}\label{fA}
f_{A}
=1-\cos^{2N_{A}}(\eps a_{A}/2)
=\frac{\eps^{2}\textcolor{black}{(a_{A})^{2}}N_{A}}{4}+\O(\eps^{4})\quad\text{as }\eps\to0.
\end{align}
To derive \eqref{fA}, note that the curved surface area of a binding site is $2\pi R_{A}^{2}(1-\cos(\eps a_{A}))$. Therefore, if $(\theta_{A}^{i},\varphi_{A}^{i})$ denotes the center of the $i$th binding site, then
\begin{align*}
f_{A}
&=1-\frac{1}{4\pi }\int_{0}^{2\pi}\int_{0}^{\pi}\P\Big((\theta,\varphi)\notin\cup_{i=1}^{N_{A}}\Gamma(\theta_{A}^{i},\varphi_{A}^{i},\eps a_{A})\Big)\sin\theta\,\dd\theta\,\dd\varphi\\
&=1-\frac{1}{4\pi }\int_{0}^{2\pi}\int_{0}^{\pi}\P\Big((\theta,\varphi)\notin\Gamma(\theta_{A}^{1},\varphi_{A}^{1},\eps a_{A})\Big)^{N_{A}}\sin\theta\,\dd\theta\,\dd\varphi\\
&=1-\frac{1}{4\pi }\int_{0}^{2\pi}\int_{0}^{\pi}\Big(1-\frac{2\pi R_{A}^{2}(1-\cos(\eps a_{A}))}{4\pi R_{A}^{2}}\Big)^{N_{A}}\sin\theta\,\dd\theta\,\dd\varphi\\
&=1-\Big[1-\tfrac{1}{2}\big(1-\cos(\eps a_{A})\big)\Big]^{N_{A}}
=1-\cos^{2N_{A}}(\eps a_{A}/2).
\end{align*}
{\black Note that the exact formula $f_{A}=1-\cos^{2N_{A}}(\eps a_{A}/2)$ in \eqref{fA} is valid for any polar angle $\eps a_{A}$ with $\eps a_{A}\in(0,\pi/2)$ and any integer $N_{A}\ge1$. In particular, one can take $N_{A}\to\infty$ with a fixed $\eps a_{A}\in(0,\pi/2)$ (which means that the binding sites necessarily overlap and cover the entire sphere) and obtain the desired result $f_{A}\to1$.}

In view of \eqref{keffr}, it remains to determine the ratio $r_{A}:=k_{\alpha}/\ks'$. In order for \eqref{keffr} to agree with the recent asymptotic results of \cite{lawley2019bp} in \eqref{recent} as $\eps\to0$, we need
\begin{align}\label{agree}
\keff\sim\eps a_{A}\lambda_{A}N_{A}/\pi\quad\text{as }\eps\to0.
\end{align}
Using \eqref{keffr}, it therefore must be the case that
\begin{align*}
r_{A}
\sim \frac{1}{\pi}\lambda_{A}N_{A}a_{A}\eps
\quad\text{as }\eps\to0,
\end{align*}
or equivalently,
\begin{align}\label{eqv}
r_{A}
\sim\frac{2}{\pi}\lambda_{A}\sqrt{N_{A}f_{A}}
\quad\text{as }f_{A}\to0.
\end{align}
Since we are interested in the case that $f_{A}\ll1$, we simply set
\begin{align}\label{set}
r_{A}
=\frac{2}{\pi}\lambda_{A}\sqrt{N_{A}f_{A}},
\end{align}
so that \eqref{keffr} becomes
\begin{align}\label{keffrn}
\keff
=\frac{f_{A}+\frac{2}{\pi}\lambda_{A}\sqrt{N_{A}f_{A}}}{1+\frac{2}{\pi}\lambda_{A}\sqrt{N_{A}f_{A}}}\ks.
\end{align}
We note that if we expand the numerator and denominator of \eqref{keffrn}, then we recover the formula obtained in \cite{lawley2019bp},
\begin{align*}
\keff
=\frac{\lambda_{A}N_{A}a_{A}\eps/\pi+\O(\eps^{2})}{1+\lambda_{A}N_{A}a_{A}\eps/\pi+\O(\eps^{2})}\ks
=\frac{\lambda_{A}N_{A}a_{A}\eps}{\pi+\lambda_{A}N_{A}a_{A}\eps}\ks+\O(\eps^{2})
=\overline{k_{A}}+\O(\eps^{2}).
\end{align*}

\subsection{Two heterogeneous molecules}

Now we consider the case of two heterogeneous molecules (i.e.\ $f_{A}\in(0,1)$ and $f_{B}\in(0,1)$). If we set $r_{A}$ as in \eqref{set} and analogously set $r_{B}$, then we obtain the following explicit bimolecular binding rate formula,
\begin{align}\label{kint}
\overline{\keff}
:=\frac{f_{A}f_{B}}{\Lambda_{A}\Lambda_{B}+\psi}\ks,\quad\text{with }r_{A}=\frac{2}{\pi}\lambda_{A}\sqrt{N_{A}f_{A}},\; r_{B}=\frac{2}{\pi}\lambda_{B}\sqrt{N_{B}f_{B}},
\end{align}
and $\Lambda_{A},\Lambda_{B},\psi$ as in \eqref{solc} and $f_{A},f_{B}$ as in \eqref{fA}.

It is straightforward to check that \eqref{kint} reduces to \eqref{keffrn} if we take $N_{B}\to\infty$ and/or $\lambda_{B}\to\infty$ (and of course the analogous statement holds if $N_{A}\to\infty$ and/or $\lambda_{A}\to\infty$). That is, \eqref{kint} has the correct limiting behavior if one or more of the four parameters, $N_{A},\lambda_{A},N_{B},\lambda_{B}$, is taken to infinity.

Does \eqref{kint} have the correct limiting behavior as $\eps\to0$? Expanding \eqref{kint} yields
\begin{align}\label{qcs}
\overline{\keff}
&=\eps^{3}N_{A}N_{B}\chiqc\ks+\O(\eps)^{4},\quad\text{as }\eps\to0,
\end{align}
where
\begin{align}\label{chia2}
\chiqc(\lambda_{A},\lambda_{B},a_{A},a_{B})
&:=\frac{a_{A}a_{B}(a_{A}\lambda_{B}+a_{B}\lambda_{A})}{4\pi}.
\end{align}
Comparing \eqref{qcs} with the behavior derived in \eqref{k0calc}, it follows that \eqref{kint} has the correct behavior as $\eps\to0$ if and only if $\chiqc=\chi$.

Using the kinetic Monte Carlo method developed in section~\ref{calculatingc0}, we find that $\chiqc\neq\chi$. However, it turns out that $\chiqc$ and $\chi$ are fairly close, $\chiqc\approx\chi$. Indeed, we plot the relative error between $\chiqc$ and $\chi$ in Figure~\ref{figchi} and find that this error is less than $16\%$ for $(D_{A},D_{B})=(R^{2}\Deff_{A}/\Dtr,R^{2}\Deff_{B}/\Dtr)\in[10^{-2},10]^{2}$ and $a_{A}=a_{B}=1$ (similar errors were found for other choices of $a_{A}$ and $a_{B}$).

While the binding rate $\overline{\keff}$ in \eqref{kint} is explicit, the formula is fairly complicated. A simpler formula that agrees with $\overline{\keff}$ in \eqref{kint} quite well is
\begin{align}\label{kints}
\overline{k_{0}}
=\frac{\eps^{3}N_{A}N_{B}}{\chi^{-1}+\eps^{2}\pi(N_{A}/(a_{B}\lambda_{B})+N_{B}/(\lambda_{A}a_{A}))+\eps^{3}N_{A}N_{B}}\ks\in(0,\ks).
\end{align}
Note that since $\chi\approx\chiqc$, for simplicity we could replace $\chi$ by $\chiqc$ in the definition of $\overline{k_{0}}$ and obtain similar results. We compare the approximations \eqref{kint} and \eqref{kints} to stochastic simulations of the full binding model in the next section.

\section{Numerical validation} \label{numerical}

In this section, we present results from two simulation methods to verify our results numerically.

\subsection{Zero rotational and surface diffusion}\label{zero}

In this subsection, we verify our asymptotic formula \eqref{k0calc} for the bimolecular binding rate in the case that $D_{A}=D_{B}\to0$. First, using the kinetic Monte Carlo method of section~\ref{calculatingc0}, we find that 
\begin{align*}
\chi(D_{A},D_{B})\approx0.1459\quad\text{for }D_{A}=D_{B}=10^{-4}.
\end{align*}
from $M=10^{10}$ trials (we take $a_{A}=a_{B}=1$). Further, the probability that $\chi(D_{A},D_{B})\in[0.1458,0.1460]$ for $D_{A}=D_{B}=10^{-4}$ is approximately $0.95$ (this follows by using the method described in section~\ref{accuracy}). Hence, the asymptotic formula \eqref{k0calc} becomes
\begin{align}\label{k0limit}
k_{0}\approx\eps^{3}N_{A}N_{B}(0.1459)\ks,\quad \eps\ll1,\,D_{A}=D_{B}\ll1. 
\end{align}

To verify \eqref{k0limit}, notice that if $D_{A}=D_{B}=0$, then the spherical caps are immobile. In particular, the problem becomes equivalent to a single point particle diffusing exterior to a 3D sphere of radius $R$ which can only be absorbed at the sphere if it reaches a pair of overlapping $A$ and $B$ spherical caps. Since the caps are immobile, if they are not initially overlapping, then they will never overlap and the particle is certain to never reach their intersection.

For simplicity, consider the case that $N_{A}=N_{B}=a_{A}=a_{B}=1$. Notice that the $B$ cap will overlap with the $A$ cap if any only if the center of the $B$ cap lies in a spherical cap of polar angle $2\eps$ centered at the $A$ cap. Assuming the caps are placed independently and uniformly on the sphere and noting that the curved surface area of a cap with polar angle $2\eps$ is $2\pi R^{2}(1-\cos(2\eps))$, the probability that the caps overlap is
\begin{align}\label{poverlap}
\frac{2\pi R^{2}(1-\cos(2\eps))}{4\pi R^{2}}
=\frac{1}{2}(1-\cos(2\eps))=\eps^{2}+\O(\eps^{4}),\quad\text{as }\eps\to0.
\end{align}

We now calculate the probability that the particle reaches the intersection of the caps, conditioned on the event, $E(s)$, that the distance (the curved geodesic distance on the sphere) between the centers of the caps is $sR\eps\ge0$, where $s\in[0,2)$. The key point is that since the caps are immobile, this problem falls into the class of problems analyzed in \cite{lindsay2017}.

Let $\P_{r}$ denote the probability measure conditioned on an initial particle radius $X(0)=r\ge R$ and an independent and uniform distribution of the initial angles $(\Theta_{0}(0),\Phi_{0}(0))$, $(\Theta_{A}^{i}(0),\Phi_{A}^{i}(0))$ for $i\in\{1,\dots,N_{A}\}$, and $(\Theta_{B}^{j}(0),\Phi_{B}^{j}(0))$ for $j\in\{1,\dots,N_{B}\}$. Recall the definition of $\tau$ in \eqref{tau}, and thus $\tau<\infty$ is the event that the particle eventually reaches the intersection of the $A$ and $B$ cap. It follows immediately from the leading order term in (3.37a) in Principal Result 3.1 in \cite{lindsay2017} that
\begin{align}\label{standard}
\P_{R_{0}}(\tau<\infty\,|\,E(s))
\sim\frac{\eps c(s)}{2R_{0}},\quad\text{as }\eps\to0,
\end{align}
where $c(s)$ is the electrostatic capacitance of the magnified ``lens'' formed by the intersection of the two spherical caps. That is, suppose $w(x,y,z;s)$ is harmonic in upper-half space,
\begin{align*}
\big(\partial_{xx}+\partial_{yy}+\partial_{zz}\big)w
&=0,\quad z>0,
\end{align*}
with mixed boundary conditions at $z=0$,
\begin{align*}
w
&=1,\quad z=0,\, (x-s/2)^{2}+y^{2}<1,\, (x+s/2)^{2}+y^{2}<1,\\
\partial_{z}w
&=0,\quad z=0,\,\text{otherwise}.
\end{align*}
Then $c(s)$ is such that
\begin{align*}
w(x,y,z;s)\sim \frac{c(s)}{\sqrt{x^{2}+y^{2}+z^{2}}},\quad\text{as }\sqrt{x^{2}+y^{2}+z^{2}}\to\infty.
\end{align*}

Now, it is straightforward to check that the probability density that the caps overlap with separation $sR\eps\in[0,2R\eps)$ given that they overlap is
\begin{align*}
\rho(s)
=\frac{s}{2}.
\end{align*}
Therefore, by conditioning on the value of the overlap distance $s\in[0,2)$ and using \eqref{poverlap} and \eqref{standard}, we obtain
\begin{align}\label{rre}
\P_{R_{0}}(\tau<\infty)
\sim\eps^{2}\int_{0}^{2}\P(R_{0}\,|\,E(s))\rho(s)\,\dd s
\sim\frac{\eps^{3}}{4R_{0}}\int_{0}^{2}c(s) s\,\dd s,\quad\text{as }\eps\to0.
\end{align}

Combining \eqref{rre} with \eqref{k0C} and \eqref{pffc}, in order to verify \eqref{k0limit} we want to show that
\begin{align}\label{want}
\frac{1}{4}\int_{0}^{2}c(s) s\,\dd s
\approx0.1459.
\end{align}
We do not have an analytic formula for $c(s)$ (except in the case $s=0$). However, we can apply the kinetic Monte Carlo method of \cite{bernoff2018b} to calculate $c(s)$ for a range of values of $s\in[0,2)$ in order to numerically compute the integral in \eqref{rre}. Taking a uniform grid of $400$ values of $s\in[0,2)$ and computing each $c(s)$ value with $10^{7}$ simulations (with an outer ``escape'' radius of $10^{5}$) yields
\begin{align*}
\int_{0}^{2}c(s) s\,\dd s
\approx0.5806.
\end{align*}
Using this numerical value for $\int_{0}^{2}c(s) s\,\dd s$, we obtain
\begin{align*}
\frac{\frac{1}{4}\int_{0}^{2}c(s) s\,\dd s}{0.1459}=0.995\approx1,
\end{align*}
which confirms \eqref{want}.

\subsection{Monte Carlo simulations of full process}

\begin{figure}[t!]
  \centering
    \includegraphics[width=1\textwidth]{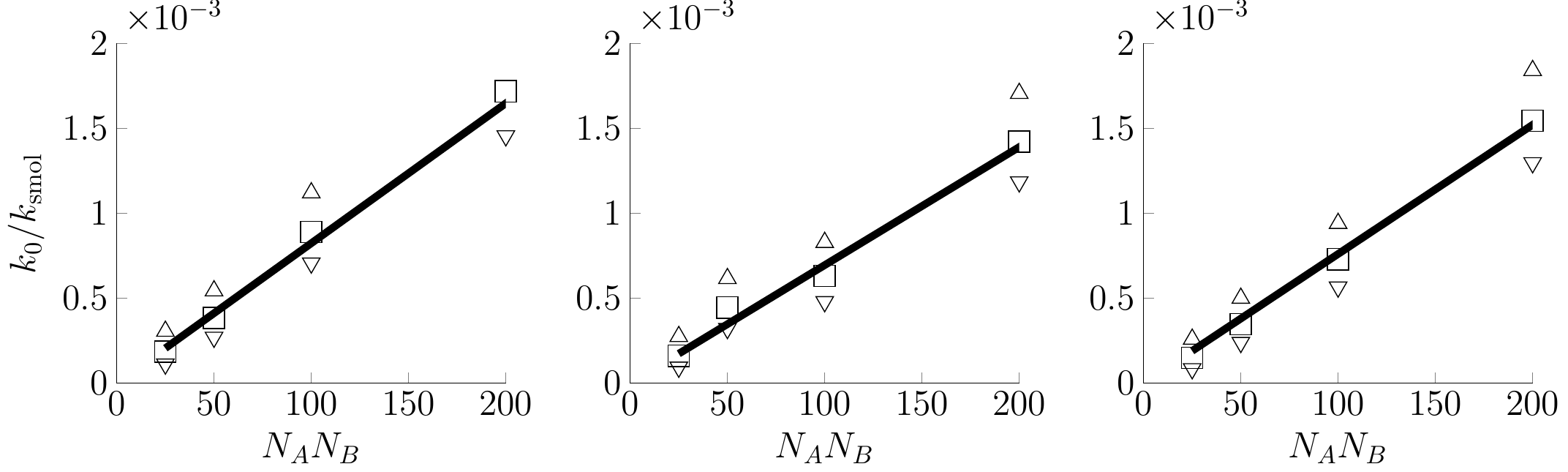}
 \caption{The asymptotic behavior of $k_{0}$ in \eqref{k0calc} as $\eps\to0$ as a function of the product $N_{A}N_{B}$ for different values of the diffusivities. In all 3 plots, the solid lines are $\eps^{3}N_{A}N_{B}\chi$ in \eqref{k0calc} where $\chi$ is computed using the kinetic Monte Carlo method of section~\ref{calculatingc0}, the squares are results from $M=10^{5}$ Monte Carlo simulations of the full process, and the triangles denote the $95\%$ confidence intervals for the simulation data using the method of section~\ref{accuracy} (all lines and data are normalized by $\ks$). In the left plot, $\Drot_A=0$, $\Drot_B=\frac{1}{2}$, $\Dsurf_A=1$, $\Dsurf_B=0$. In the middle plot, $\Drot_A=\Drot_B=\Dsurf_A=\Dsurf_B=\frac{1}{2}$. In the right plot, $\Drot_A=\Drot_B=\frac{1}{2}$, $\Dsurf_A=\Dsurf_B=0$.
 }
 \label{figmain}
\end{figure}

In this subsection, we compare the asymptotic formula for $k_{0}$ in \eqref{k0calc} and the approximations $\overline{\keff}$ and $\overline{k_{0}}$ (equations \eqref{kint} and \eqref{kints}) to Monte Carlo simulations of the full process in section~\ref{stochastic}.  Before describing our stochastic simulation method in more detail, we first outline the main points and give the results.

As shown in section~\ref{math}, the problem is equivalent to (i) a set of $N_{A}$-many $A$ spherical caps and a set of $N_{B}$-many $B$ spherical caps that each move on the surface of a single sphere with radius $R$ and (ii) a point particle that diffuses exterior to this sphere and is absorbed at the sphere if and only if it hits the intersection of an $A$ spherical cap with a $B$ spherical cap (otherwise it reflects from the sphere). The motion of the $A$ spherical caps is governed by their individual surface diffusions (with surface diffusivity $\Dsurf_{A}$) and the rotational diffusion of the $A$ molecule (with rotational diffusivity $\Drot_{A}$), and similarly for $B$ spherical caps.

We thus simulate the path of a single particle with diffusivity $\Dtr$ in $\R^{3}$ exterior to a sphere of radius $R>0$ and the paths of the diffusing caps on the surface of the sphere until the particle either reaches the intersection of the $A$ and $B$ caps or reaches some large outer radius $R_{\infty}\in(R,\infty)$. After repeating this $M\gg1$ times, we calculate the proportion of particles that reach $R_{\infty}$. A certain modification of this proportion then yields an approximation to the probability $p$ in \eqref{pdef} that the particle never reaches the intersection of the $A$ and $B$ caps, which then yields an approximation to $k_{0}$ via \eqref{pffc}-\eqref{k0C}.

\begin{figure}[t!]
  \centering
    \includegraphics[width=1\textwidth]{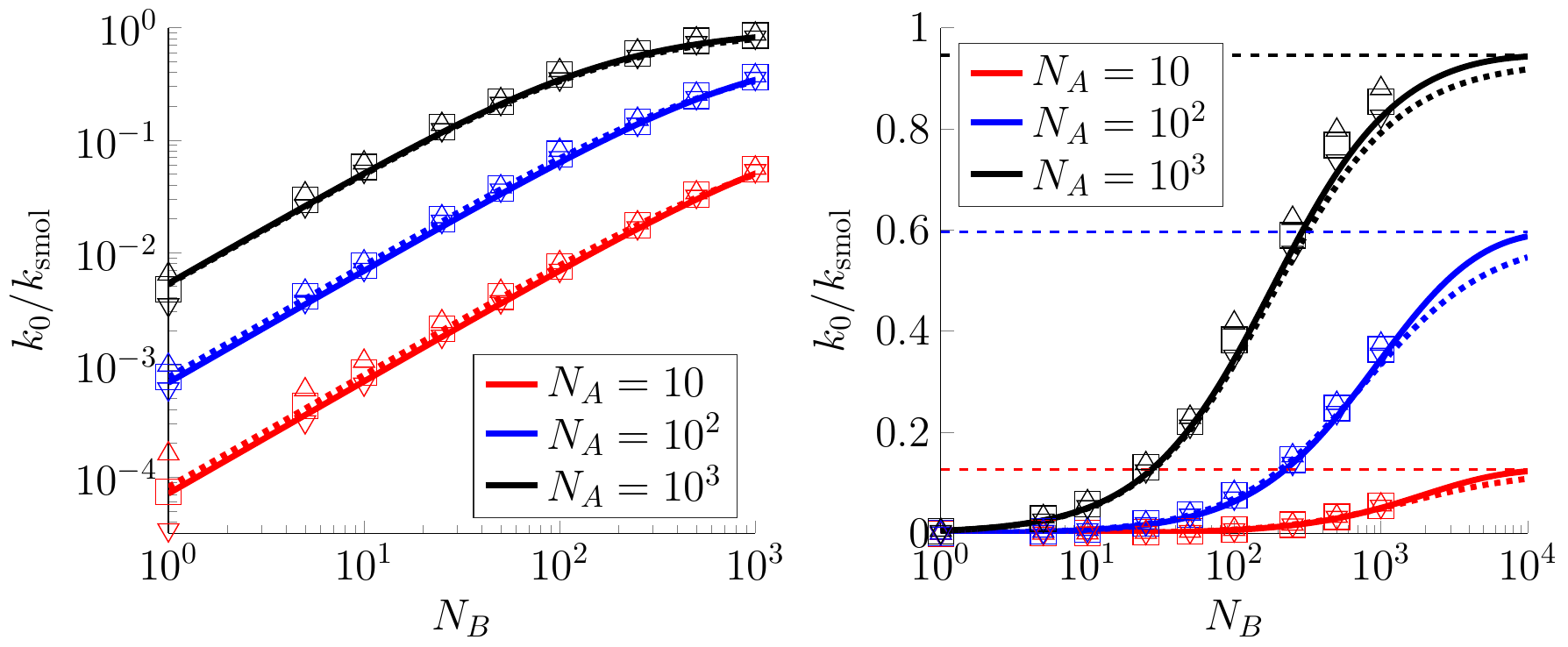}
 \caption{The approximations $\overline{\keff}$ and $\overline{k_{0}}$ to the bimolecular binding rate $k_{0}$ as a function of $N_{B}$ for $N_{A}\in\{10,10^{2},10^{3}\}$. In both plots, the solid curves are $\overline{\keff}$ in \eqref{kint}, the dotted curves are $\overline{k_{0}}$ in \eqref{kints}, the squares are results from Monte Carlo simulations of the full process, and the triangles denote the $95\%$ confidence intervals for the simulation data using the method of section~\ref{accuracy} (all curves and data are normalized by $\ks$). The left and right plots have identical data, but the left plot has a logarithmic vertical axis for better visualization for small values of $N_{B}$. In the right plot, the horizontal dashed lines give $\lim_{N_{B}\to\infty}\overline{\keff}$. In both plots, we take $\Drot_{A}=\Drot_{B}=1$ and $\Dsurf_{A}=\Dsurf_{B}=0$.
 }
 \label{figNs}
\end{figure}

Figures~\ref{figmain}-\ref{figNs} show very good agreement between our theoretical results and these stochastic simulations of the full process. Figure~\ref{figmain} plots the asymptotic formula for $k_{0}$ in \eqref{k0calc} for different values of the diffusivities $\Drot_{A}$, $\Drot_{B}$, $\Dsurf_{A}$, and $\Dsurf_{B}$ as a function of the product $N_{A}N_{B}\in[25,200]$, where $\chi$ is computed from the kinetic Monte Carlo simulations of section~\ref{calculatingc0}.

Since the formula $k_{0}\sim\eps^{3}N_{A}N_{B}\chi\ks$ of \eqref{k0calc} clearly breaks down for fixed $\eps$ and large $N_{A}$ or $N_{B}$, Figure~\ref{figNs} plots the approximations $\overline{\keff}$ and $\overline{k_{0}}$ in \eqref{kint} and \eqref{kints} for larger values of $N_{A}$ and $N_{B}$. From this figure, we see that $\overline{\keff}\approx\overline{k_{0}}$ agrees well with stochastic simulations. The largest errors tend to occur when $N_{A}=10^{3}$ (or $N_{B}=10^{3}$), in which case the surface area fraction covered by binding sites is $f_{A}\approx N\eps^{2}/4=\frac{1}{4}$ (or $f_{B}\approx\frac{1}{4}$). This error is expected, since our approximations were made assuming the surface area fraction is small. We note that the data points in Figure~\ref{figNs} were computed from either $M=10^{3}$, $M=10^{4}$, or $M=10^{5}$ simulations of the full process, depending on the values of $N_{A}$ and $N_{B}$. In particular, we used a smaller number of simulations for large values of $N_{A}$ and/or $N_{B}$ as such simulations are very computationally expensive (fortunately, larger values of $N_{A}$ and/or $N_{B}$ yield a higher binding probability, so less simulations are needed to get a precise estimate of $k_{0}$). We further note that we take $\Drot_{A}=\Drot_{B}=1$ and $\Dsurf_{A}=\Dsurf_{B}=0$ in Figure~\ref{figNs} to decrease computational cost.

We now describe our simulation method, which is similar to the method used in \cite{lawley2019bp,lawley2019dtmfpt,lawley2019fpk}. Initially, we place the particle at radius $R_{0}\in(R,R_{\infty})$ and randomly distribute the caps uniformly on the sphere. The diffusion of the particle (with diffusivity $\Dtr$), the surface diffusion of the $A$ caps (with diffusivity $\Dsurf_{A}$), and the surface diffusion of the $B$ caps (with diffusivity $\Dsurf_{B}$) are simulated with the Euler-Maruyama method \cite{kloeden1992}. To increase computational efficiency, we use either a large time step (denoted $\Delta t_{\text{big}}$) or a small time step (denoted $\Delta t_{\text{small}}$), depending on the distance between the particle and the nearest $A$ and $B$ caps. 

To implement rotational diffusion, at each time step, all the $A$ caps undergo the same random rotations about the three Cartesian coordinate axes (and similarly for the $B$ caps). More precisely, if we define the rotation matrices,
\begin{align*}
R_{x}(\omega)
&:=\begin{pmatrix}
1 & 0 & 0\\
0 & \cos\omega & -\sin\omega\\
0 & \sin\omega & \cos\omega
\end{pmatrix},
\quad
R_{y}(\omega)
:=\begin{pmatrix}
\cos\omega & 0 & \sin\omega\\
0 & 1 & 0\\
-\sin\omega & 0 & \cos\omega
\end{pmatrix}
,\\
R_{z}(\omega)
&:=\begin{pmatrix}
\cos\omega & -\sin\omega & 0\\
\sin\omega & \cos\omega & 0\\
0 & 0 & 1
\end{pmatrix},
\end{align*}
and let $\{(x_{A}^{i},y_{A}^{i},z_{A}^{i})\}_{i=1}^{N_{A}}$ denote the Cartesian coordinates of the centers of the $A$ caps at the start of a time step of size $\Delta t>0$, then the centers of these caps at the end of the time step are
\begin{align*}
R_{x}(\omega_{1})R_{y}(\omega_{2})R_{z}(\omega_{3})\begin{pmatrix}
x_{A}^{i}\\ y_{A}^{i}\\ z_{A}^{i}
\end{pmatrix}\in\R^{3},\quad i\in\{1,\dots,N_{A}\},
\end{align*}
where $\omega_{1},\omega_{2},\omega_{3}$ are 3 independent realizations of Gaussian random variables with mean zero and variance $2\Drot_{A}\Delta t>0$. Importantly, the random variables $\omega_{1},\omega_{2},\omega_{3}$ do not depend on the index $i\in\{1,\dots,N_{A}\}$, which means that the random rotation is common to all of the $A$ caps. The $B$ caps are rotated analogously, though of course the $B$ caps are rotated independently from the $A$ caps.

In all simulations, we take $\eps=10^{-1.5}\approx0.03$, $\Dtr=R=a_{A}=a_{B}=1$, $R_{0}=1.1$, $R_{\infty}=10$, $\Delta t_{\text{big}}=10^{-3}$, $\Delta t_{\text{small}}=10^{-8}$, and $M\in[10^{3},10^{5}]$ trials.

We now describe more precisely how we estimate $k_{0}$ from the simulation data. Let $\P_{r}$ denotes the probability measure conditioned on an initial particle radius $X(0)=r\ge R$ and an independent and uniform distribution of the initial angles $(\Theta_{0}(0),\Phi_{0}(0))$, $(\Theta_{A}^{i}(0),\Phi_{A}^{i}(0))$ for $i\in\{1,\dots,N_{A}\}$, and $(\Theta_{B}^{j}(0),\Phi_{B}^{j}(0))$ for $j\in\{1,\dots,N_{B}\}$. For a set of $M\gg1$ trials with fraction $q\in[0,1]$ that reach radius $R_{\infty}$, we obtain the approximation
\begin{align}\label{q}
q\approx\P_{R_{0}}(\tau_{R_{\infty}}<\tau),
\end{align}
where $\tau\ge0$ is the first time the particle reaches the intersection of $A$ and $B$ caps (defined in \eqref{tau}), and $\tau_{R_{\infty}}\ge0$ is the first time the particle reaches radius $R_{\infty}$,
\begin{align*}
\tau_{R_{\infty}}:=\inf\{t>0:X(t)=R_{\infty}\}.
\end{align*}
Now, it follows immediately from integrating \eqref{ppde} and using \eqref{pffc} that
\begin{align*}
1-\tfrac{C}{R_{0}}
=\P_{R_{0}}(\tau=\infty),\quad R_{0}>R.
\end{align*}
Therefore, to find an approximation for $\P_{R_{0}}(\tau=\infty)$, we follow \cite{Northrup1988, Northrup1986, lawley2019bp} to obtain
\begin{align}
1-\tfrac{C}{R_{0}}
&=\P_{R_{0}}(\tau=\infty\,|\,\tau>\tau_{R_{\infty}})\P_{R_{0}}(\tau>\tau_{R_{\infty}})\nonumber\\
&\approx\P_{R_{\infty}}(\tau=\infty)\P_{R_{0}}(\tau>\tau_{R_{\infty}})\label{app}
=\big(1-\tfrac{C}{R_{\infty}}\big)\P_{R_{0}}(\tau>\tau_{R_{\infty}}).
\end{align}
The error in the approximation in \eqref{app} vanishes as $R_{\infty}/R$ and/or $\min\{\Deff_{A},\Deff_{B}\}/\Dtr$ grows. If we rearrange \eqref{app} and use \eqref{q}, then we obtain the following numerical approximation to the capacitance in \eqref{pffc},
\begin{align*}
C\approx \frac{(1-q)R_{0}R_{\infty}}{R_{\infty}-qR_{0}}.
\end{align*}
Plugging this approximation for $C$ into \eqref{k0C} yields the numerical approximation to the bimolecular binding rate $k_{0}$ that is used in plotting the ratio $k_{0}/\ks$ in Figures~\ref{figmain}-\ref{figNs}. 

\section{Discussion} \label{discussion}

In this paper, we considered a generalization of the classical Smoluchowski model for bimolecular binding rates that includes the fact that pairs of molecules can bind only in certain orientations. This generalization took the form of a high-dimensional, anisotropic diffusion equation with mixed boundary conditions. We applied matched asymptotic analysis \cite{lindsay2017} to this PDE and derived the bimolecular binding rate in the limit of small binding sites. The resulting binding rate formula involves a factor, $\chi$, that we computed numerically by modifying a recent kinetic Monte Carlo algorithm \cite{bernoff2018b}. We then applied the quasi chemical approximation \cite{solc1973} to obtain (i) a formula which includes the effects of binding site competition/saturation and (ii) a simple analytical approximation for $\chi$. Our analysis thus constitutes a hybrid asymptotic-numerical approach \cite{kropinski1995, lindsay2016, wang2019}, as it relied on both asymptotic analysis and numerical computation.

In our model, both particles are ``patchy'' or ``heterogeneous,'' meaning that both particles contain localized binding sites. The limiting case of one heterogeneous molecule and one homogeneous molecule (one molecule completely covered in binding sites) is a classical and well-studied problem \cite{lindsay2017,Berezhkovskii2004,Berezhkovskii2006,Dagdug2016,Eun2017,lindsay2017,Muratov2008,Zwanzig1990,eun2020,kaye2019}, dating back to Berg and Purcell's landmark 1977 work \cite{berg1977} which yielded the rate constant,
\begin{align*}
\kbp
:=\frac{\eps a_{A}N_{A}}{\pi+\eps a_{A}N_{A}}\ks.
\end{align*}
A number of interesting works have modified Berg and Purcell's formula to account for the effects of binding site arrangement and curvature of the molecular surface \cite{lindsay2017,Berezhkovskii2004,Berezhkovskii2006,Dagdug2016,Eun2017,lindsay2017,Muratov2008,Zwanzig1990,eun2020,lawley2019bp,kaye2019}. In fact, the method of matched asymptotic analysis that we employed in the present work follows the method employed in \cite{lindsay2017}, and also similar methods in \cite{cheviakov11,cheviakov2010,coombs09,lawley2019fpk, lawley2019dtmfpt,PB3,PB4,lawley_asymptotic_2018}. These formal methods are related to the strong localized perturbation analysis pioneered in \cite{ward93,ward93b}. Note that the model of the present work is a strict generalization of the model of Berg and Purcell.

The model of two heterogeneous molecules that we analyzed in the present work was studied in the case of a single binding site on each molecule ($N_{A}=N_{B}=1$) by {\v{S}}olc and Stockmayer in 1971 \cite{solc1971}. In that work, the authors used the symmetry inherent in the single binding site model to derive an expression for the binding rate in terms of an infinite series requiring the solution of an infinite system of linear algebraic equations \cite{solc1971}. In the absence of a tractable expression for the binding rate for this single binding site model, subsequent studies have employed either the so-called closure (constant flux) approximation \cite{temkin1984, lee1987} or the quasi chemical approximation \cite{solc1973}. Though heuristic, the quasi chemical approximation was shown to be quite accurate for the case of a single, relatively large binding site on each molecule \cite{zhou1993}. In the case of a single small binding site on each molecule, the quasi chemical approximation \cite{solc1973} combined with the analysis of Berg \cite{Berg1985} predicts that the bimolecular binding rate has the following approximate asymptotic behavior \cite{zhou1993},
\begin{align}\label{berga}
k_{0}
\approx
\eps^{3}\chib\ks,\quad\text{if }\eps\ll1,
\end{align}
where
\begin{align}\label{chib}
\chib=\chib(\lambda_{A},\lambda_{B},a_{A},a_{B})
&:=\frac{a_{A}a_{B}(a_{A}\lambda_{B}+a_{B}\lambda_{A})}{8\sqrt{2}}.
\end{align}
Comparing \eqref{chib} with our formula for $\chiqc$ in \eqref{chia2} (which approximates the quantity $\chi$ determined numerically in section~\ref{calculatingc0}), we see that the only difference is the factor $1/(8\sqrt{2})\approx0.09$ in \eqref{chib} versus the factor $1/(4\pi)\approx0.08$ in \eqref{chia2}. This difference arises because \eqref{chib} relies on an approximation of a certain infinite series, whereas \eqref{chia2} depends on the asymptotic predictions of \cite{lawley2019bp} (see (17)-(18a) in \cite{Berg1985} and the discussion surrounding (61) in \cite{lawley2019bp} for more details). Hence, the results in this paper show that the heuristic prediction \eqref{berga} is quite accurate as $\eps\to0$ and extend the binding rate formula to the case of multiple binding sites.


Related work that studied the binding of spherical molecules with multiple binding sites (often called molecules with ``patches'' or simply ``patchy particles'') includes \cite{newton2015, roberts2014, klein2014}. In particular, reference \cite{newton2015} used Monte Carlo simulations to investigate the relative contributions of translational and rotational diffusion to the association of two or more patchy particles.  Reference \cite{roberts2014} studied the association of pairs of patchy particles with a few relatively large patches using lattice models and lattice-adjacent models, and reference \cite{klein2014} introduced a computational approach for studying association and dissociation of such patchy particles. {\black
In addition to models with spherical molecules, progress has recently been made in analytically studying diffusion-influenced reactions for non-spherical molecules \cite{galanti2016, traytak2018, piazza2019}.} 


{\black
In previous work, the mixed boundary conditions that result from patchy particles are often approximated by a homogeneous boundary condition through a method called boundary homogenization \cite{Berezhkovskii2004}. Specifically, one considers the Smoluchowski problem in \eqref{pdeG}-\eqref{ffc} with the absorbing boundary condition \eqref{acG} at the reaction radius $r=R$ replaced by a Robin boundary condition (also called a partially absorbing, radiation, third type, impedance, or convective condition \cite{lawley15jk, chapman2016, erban_reactive_2007, singer_partially_2008}),
\begin{align}\label{robin}
D\partial_{r}p
=\kappa p,\quad \text{for }r=R,
\end{align}
where $\kappa>0$ is the so-called trapping rate parameter (or partial reactivity). Cast in this form, the problem becomes one of choosing the homogenized trapping rate in order to approximate the heterogeneous reactivity. Solving \eqref{pdeG}-\eqref{ffc} with \eqref{robin} yields the following expression for the binding rate $k$ in terms of the trapping rate $\kappa$,
\begin{align}\label{kkappa}
k
:=\Dtr\int_{r=R}\partial_{r}p\,\dd S
=\frac{\kappa R}{\Dtr+\kappa R}\ks.
\end{align}
Solving \eqref{kkappa} for $\kappa$ yields
\begin{align}\label{kkappa2}
\kappa
=\frac{\Dtr}{R}\Big(\frac{k/\ks}{1-k/\ks}\Big).
\end{align}
Hence, while our results are given in terms of binding rates (i.e.\ $k_{0}$, $\overline{k_{0}}$, $\keff$, etc.), these formulas can be translated into the corresponding trapping rates via \eqref{kkappa2}.
}

{\black Related to the point above, we note that our model assumes that the molecules bind immediately once the binding sites are in contact. Mathematically, this assumption manifests as the absorbing Dirichlet condition in the mixed Dirichlet/Neumann boundary conditions in \eqref{pbc}. This assumption can be relaxed by replacing the Dirichlet/Neumann conditions in \eqref{pbc} by Robin/Neumann conditions of the form
\begin{alignat*}{2}
D\partial_{r}p
&=\kappa_{\text{bs}}p,\quad &&r=R,\,(\theta_{0},\varphi_{0})\in \Lambda(\theta_{\f},\varphi_{\f})\cap\Lambda(\theta_{\m},\varphi_{\m}),\\
\partial_{r}p
&=0, \quad &&r=R,\,(\theta_{0},\varphi_{0})\notin \Lambda(\theta_{\f},\varphi_{\f})\cap\Lambda(\theta_{\m},\varphi_{\m}),
\end{alignat*}
where $\kappa_{\text{bs}}\in(0,\infty)$ models a finite binding rate of binding sites that are in contact ($\kappa_{\text{bs}}$ is not to be confused with $\kappa$ in \eqref{robin}). Such mixed Robin/Neumann conditions are known to affect the leading order behavior in narrow escape problems involving small targets \cite{grebenkov2017, lawley2019imp}. Alternatively, a finite binding rate can be modeled by retaining the absorbing Dirichlet condition in \eqref{pbc} and reducing the binding site size to some ``effective'' size. This latter perspective is the one taken by Berg and Purcell \cite{berg1977}.
}


In closing, we briefly discuss our results in the context of empirical binding rates. The Smoluchowski bimolecular binding rate \eqref{ksmol} for typical proteins is roughly {\black \cite{northrup1992, zhou2010, zhou2013}}
\begin{align*}
\textcolor{black}{
\ks\approx7\times10^{9}\;\text{M}^{-1}\text{sec}^{-1}.
}
\end{align*}
This rate significantly overestimates experimentally measured rates, which is to be expected since it ignores orientational constraints in binding. Indeed, empirical rates are often in the range \cite{northrup1992}
\begin{align}\label{emprange}
\ke\in[0.5,5]\times10^{6}\;\text{M}^{-1}\text{sec}^{-1}.
\end{align}
As noted in the Introduction, it is tempting to account for the orientational constraints by simply multiplying the Smoluchowski rate by a geometric factor given by the product of the protein surface area fractions covered by binding sites, which yields the binding rate $\kg:=f_{A}f_{B}\ks$ (see \eqref{fsa}-\eqref{kgeo}). However, this simple modification yields a binding rate that is typically a few orders of magnitude smaller than experimentally measured rates. For example, it has been estimated that \cite{northrup1992}
\begin{align}\label{kgeonumber}
\kg\approx7\times10^{2}\;\text{M}^{-1}\text{sec}^{-1}.
\end{align}
Since $\ks$ overestimates $\ke$ and $\kg$ underestimates $\ke$, it is interesting to note that the binding rate, $k_{0}$, satisfies 
\begin{align*}
\kg\ll k_{0}\ll \ks,
\end{align*}
in the limit of small binding sites, $\eps\ll1$. Indeed, $\kg$ for our model is
\begin{align*}
\kg
=\Big(\frac{N_{A}\eps^{2}}{4}\Big)\Big(\frac{N_{B}\eps^{2}}{4}\Big)\ks
=\eps^{4}N_{A}N_{B}\frac{1}{16}\ks,
\end{align*}
where we have taken $a_{A}=a_{B}=1$ for simplicity. Hence, \eqref{k0calc} gives
\begin{align*}
\frac{\kg}{k_{0}}
\sim\eps\frac{\chi}{16}\ll1,\quad \text{for }\eps\ll1.
\end{align*}
Therefore, if we take the value \eqref{kgeonumber} for $\kg$ and for definiteness take $\chi=0.29$ from \eqref{chivals}, then we obtain that $k_{0}$ is in the typical empirical range \eqref{emprange} if $\eps\in[10^{-2},10^{-1}]$.


\bibliography{library}
\bibliographystyle{siam}
\end{document}